\documentclass[11pt, reqno]{amsart}

\usepackage[T1]{fontenc}
\usepackage[margin=2.5cm]{geometry}
\usepackage{amsmath, amssymb, amsthm, mathtools}
\usepackage{mathrsfs}
\usepackage{enumerate}
\usepackage{float}
\usepackage{url}
\usepackage{tikz-cd}
\usepackage[backend=biber, style=alphabetic, doi=false, isbn=false, url=false]{biblatex}
\usepackage{extarrows}
\usepackage{hyperref}

\newtheorem{theorem}{Theorem}[section]
\newtheorem{corollary}[theorem]{Corollary}
\newtheorem{lemma}[theorem]{Lemma}
\newtheorem{proposition}[theorem]{Proposition}
\newtheorem{conjecture}[theorem]{Conjecture}

\theoremstyle{definition}
\newtheorem{definition}[theorem]{Definition}
\newtheorem{remark}[theorem]{Remark}

\DeclareMathOperator{\Img}{Im}

\DeclareMathOperator{\Spec}{Spec}

\newcommand{\Fl}{\mathcal{F}l}

\DeclareMathOperator{\Gal}{Gal}

\DeclareMathOperator{\Tr}{Tr}

\DeclareMathOperator{\Rep}{Rep}
\DeclareMathOperator{\Ind}{Ind}
\DeclareMathOperator{\Sym}{Sym}
\DeclareMathOperator{\Hom}{Hom}
\DeclareMathOperator{\End}{End}

\DeclareMathOperator{\Isom}{Isom}
\DeclareMathOperator{\Res}{Res}
\DeclareMathOperator{\Lie}{Lie}

\DeclareMathOperator{\GL}{GL}

\DeclareMathOperator{\GU}{GU}

\newcommand{\Gm}{\mathbb{G}_m}

\DeclareMathOperator{\Fil}{Fil}

\newcommand{\C}{\mathbb{C}}
\newcommand{\R}{\mathbb{R}}
\newcommand{\Q}{\mathbb{Q}}

\newcommand{\Qp}{\mathbb{Q}_p}

\newcommand{\Z}{\mathbb{Z}}
\newcommand{\Zp}{\mathbb{Z}_p}
\newcommand{\Zpt}{\mathbb{Z}_{p^2}}
\newcommand{\Fpb}{\overline{\mathbb{F}}_p}
\newcommand{\Fp}{\mathbb{F}_p}
\newcommand{\Fpt}{\mathbb{F}_{p^2}}

\newcommand{\adf}{\mathbb{A}_{f}}
\newcommand{\adfp}{\mathbb{A}_{f}^p}

\newcommand{\M}{\mathcal{M}}

\DeclareMathOperator{\Sh}{Sh}
\DeclareMathOperator{\Gzip}{G-Zip}
\DeclareMathOperator{\Gzf}{G-ZipFlag}
\newcommand{\ShKtwo}{\Sh_K}

\newcommand{\FlKtwo}{\Fl_{\, K}}
\newcommand{\ShKtr}{\Sh_{K'}'}
\newcommand{\FlKtr}{\Fl_{\, K'}'}
\newcommand{\ShKp}{\Sh_{K, \Fpt}}
\DeclareMathOperator{\Hdr}{\mathcal{H}^1_{dR}}
\DeclareMathOperator{\Zcb}{\overline{\mathcal{Z}}}

\newcommand{\Filc}{\Fil_{conj}}

\DeclareMathOperator{\ha}{ha}

\title{An exotic congruence for irregular weight Picard modular forms}
\author{Calle Sönne}

\begin{document}

\begin{abstract}
  We construct a new type of congruence from irregular weight mod $p$ Picard modular forms (for $p$
  inert) to ones
  of regular weight, analogous to the one obtained by Kansal--Levin--Savitt \cite{KLS} in the Hilbert setting.
  The congruence is obtained by using a geometric Jacquet--Langlands phenomenon, as first established by Vollaard \cite{vollaard}, to transfer irregular weight mod $p$ Picard modular forms from the Picard modular surface to a Shimura set.
  We call this an \textit{exotic} congruence, since it relates mod $p$ automorphic forms
  associated to different groups.
  This congruence will be used in upcoming work of the author proving mod $p$ local-global compatibility
  for irregular weight Picard modular forms (in the $p$ inert case).
\end{abstract}

\maketitle
\setcounter{tocdepth}{1}
\tableofcontents

\section{Introduction}
In this paper, we establish a new type of congruence from irregular
weight Picard modular forms to Picard modular forms of regular weight.
In upcoming work of the author, this will be used to prove \textit{mod $p$ local-global compatibility}
for irregular weight Picard modular forms in the $p$ inert case,
analogous to the upcoming work of Kansal--Levin--Savitt \cite{KLS} in the Hilbert setting (for $p$ unramified).
Before describing the specific setting of the paper in more detail, let
us give some background on the problem in the simpler case of the modular curve.
For the duration of the paper, we fix a prime $p > 2$.

\subsection{Mod $p$ modular forms}
Let $X / \Spec \Zp$ denote the compactified modular curve (of some level prime to $p$).
Let $\omega$ denote the Hodge bundle associated to the universal family $\mathcal{E} / X$.
Then, as first introduced by Katz \cite{katz}, we can consider mod $p$ modular forms of weight $k \in \mathbb{Z}_{\ge 0}$ as elements of $H^0(X_{\Fpb}, \omega^{\otimes k})$.

A modular form of \textit{irregular weight} is just a weight one modular form.
These are harder to study than their regular weight (i.e. $k \ge 2$) counterparts.
One of the main tools we have for studying weight one modular form is by using
\textit{congruences} to regular weight modular forms. Here,
we say that two mod $p$ eigenforms (of potentially different weights) are congruent
if their mod $p$ Hecke eigensystems are the same at almost all primes.
In particular, the Galois representations associated to two congruent mod $p$ modular forms
will be the same.

\subsection{Mod $p$ local-global compatibility for weight one modular forms}

The following theorem, which follows from a simplified version of \cite[Theorem 3.11]{calegari-geraghty}, can
be seen as an instance of mod $p$ local-global compatibility for weight one modular forms.
\begin{theorem}[\cite{calegari-geraghty}] \label{thm:mod-p-LGC}
  Given a weight one modular form $f \in H^0(X_{\Fpb}, \omega)$, let $\rho_f$ denote the associated Galois representation. Then $\rho_f|_{G_{\Qp}}$
  is unramified.
\end{theorem}
Let us give some remarks on why this theorem (and the corresponding generalizations for other groups) is
interesting.
\begin{remark}
  \begin{enumerate}
    \item Theorems of this type are important for extending modularity lifting theorems beyond the regular weight case.
      Namely, the correct Galois deformation rings to consider will specify local conditions at $p$.
      In the case of weight one modular forms, the local condition is that the Galois representation
      is unramified at $p$.
      Hence, Theorem \ref{thm:mod-p-LGC} (or rather the full statement of \cite[Theorem 3.11]{calegari-geraghty}) is required for constructing the morphism ``$R \rightarrow \mathbb{T}$''.
    \item This question is closely related to extensions of the weight part of Serre's conjecture
      beyond the regular weight case. For modular forms, this was first considered (and proved) in
      \cite{Edix92}. Recently, Diamond--Sasaki \cite{DS} formulated an analogous extension in the setting
      of
      Hilbert modular forms.
    \item The most interesting case of this theorem is when $f$ does not admit a lift to
      characteristic $0$. It is crucial to prove the theorem also in these cases if one wants to use the method of
      \cite{calegari-geraghty} to obtain modularity lifting theorems.
  \end{enumerate}
\end{remark}

\subsection{The work of Kansal--Levin--Savitt}
In their upcoming work \cite{KLS}, the authors solve the analogous question in the Hilbert setting (assuming $p$ is unramified in the totally real field).
This means that the Shimura varieties they work with are Hilbert modular varieties,
and the automorphic forms are Hilbert modular forms.
The irregular weight Hilbert modular forms are the \textit{partial weight one} Hilbert modular forms.

The approach goes as follows. The first step is to establish multiple congruences from
irregular weight Hilbert modular forms to ones of regular weight.
Once they have these congruences, they then prove the conjecture by an argument purely on the Galois side.
In order to have enough congruences in all cases, they establish a new, third, type of congruence beyond the
well-known congruences obtained from \textit{partial Hasse invariants} and \textit{theta operators}.
This new congruence is obtained by transferring the Hilbert modular form to a quaternionic Shimura variety of smaller dimension, using a \textit{geometric Jacquet--Langlands} phenomenon.
This produces a quaternionic Hilbert modular form,
with the same mod $p$ Hecke eigensystem (away from $p$), and of regular weight.
\begin{remark}
  Let us mention some other relevant work in the Hilbert setting.
  First, Wiersema \cite{wiersema}
  proved mod $p$ local-global compatibility for Hilbert modular forms in 
  the cases when this third type of congruence is not needed. This was inspired by earlier work of
  Diamond--Sasaki \cite{DS} in the inert quadratic case.
  There is also the work of Dimitrov--Wiese \cite{dimitrov-wiese}, who proved mod $p$ local-global compatibility
  in the case of parallel weight one.
\end{remark}

\subsection{Mod $p$ Picard modular forms}
The goal of this paper is to establish an analogue of this new type of congruence
for mod $p$ Picard modular forms of irregular weight.
Since this will really be a congruence between Picard modular forms on different
Shimura varieties, we propose to call this an \textit{exotic} congruence.
Let us now briefly introduce the Shimura variety of interest, and the spaces of mod $p$ Picard modular forms.
We refer to the main text for more details.
\subsubsection{The Shimura variety} 
Fix an imaginary quadratic field $E / \Q$ such that $p$ is inert in $E$.
Our Shimura variety $\Sh_K / \Spec \Fpt$ (of some level $K$, but we ignore this during the introduction) will be the moduli space of abelian threefolds equipped with a prime-to-$p$ quasi-polarization, an $\mathcal{O}_E$-action and satisfying the Kottwiz determinant condition of signature $(2, 1)$.
Let $\mathcal{A} / \Sh_K$ denote the universal abelian variety.

We will also be considering the so-called \textit{flag Shimura varieties}. These
are useful as they allow us to describe mod $p$ Picard modular forms as sections
of line bundles, whereas on the Shimura variety they can only be described as sections of vector bundles in general.
Let us remark that for the modular curve, the ``flag modular curve'' is just the modular curve itself.
We let $\Fl_K / \Spec \Fpt$ denote the corresponding flag Shimura variety. It comes equipped with a morphism $\pi :
\Fl_K \rightarrow \Sh_K$, whose fibers are pointwise isomorphic to $\mathbb{P}^1$.

\subsubsection{The Shimura set}
We will also consider a $0$-dimensional \textit{Shimura set} $\ShKtr / \Spec \Fpt$ (of some level
$K'$, which we also ignore for now).
It has the same moduli interpretation as $\Sh_K$, except we take the signature to be $(3, 0)$.
Let $\FlKtr$ denote the associated flag Shimura variety.
There is again a map $\pi' : \FlKtr \rightarrow \ShKtr$, but this time it has relative dimension $3$, and the fibers are given by flag varieties.

\subsubsection{Mod $p$ Picard modular forms}
Choose a labeling $\Hom(\mathcal{O}_{E_p}, \Zpt) = \{ \tau_0, \tau_1 \}$.
Let $\omega_{\mathcal{A}, \tau_0}$ denote the $\tau_0$-isotypic component of the Hodge bundle of $\mathcal{A}/\Sh_K$.
The signature condition implies that this is of rank $2$. We can then define the \textit{automorphic
vector bundles} $\omega(k, l) := \Sym^{k} \omega_{\mathcal{A}, \tau_0} \otimes {\det}^{l}
\omega_{\mathcal{A}, \tau_0}$ on
$\ShKtwo$ for $k, l \in \Z$ and $k \ge 0$ (for simplicity we only use two parameters here, see \S\ref{sect:aut-vb} for a more general definition).
A mod $p$ Picard modular form of weight $\lambda = (k, l)$ is then given by an element of $H^0(\Sh_{K, \Fpb}, \omega(\lambda))$.
Note that we do not need to consider compactifications here by Koecher's principle
\cite{koecher}.
We can also define a line bundle $\mathcal{L}(\lambda)$ on $\FlKtwo$
such that $H^0(\Fl_{K, \Fpb}, \mathcal{L}(\lambda)) = H^0(\Sh_{K, \Fpb}, \omega(\lambda))$.

In our setting, a weight $(k, l)$ will be \textit{irregular} if $l = 2$.
The point is that these are the weights for which, in characteristic $0$,
the automorphic representations contributing to the cuspidal cohomology will
be holomorphic limits of discrete series at $\infty$
(see e.g. \cite[§2]{goldring} \cite[§3]{harris} \cite{BHR}).

We can also consider Picard modular forms on the Shimura set $\ShKtr$.
Given any tuple $\lambda = (a, b) \in
\Z^{\oplus 2}_{\ge 0}$ one can define an automorphic vector bundle $\omega'(\lambda)$ on $\ShKtr$.
Note that we only use $2$ parameters to describe the weight for simplicity, and throughout the
introduction we will be ignoring issues with determinant and multiplier twists.
The global sections of $\omega'(\lambda)$ can be identified with spaces of \textit{algebraic automorphic forms} in the sense of Gross \cite{gross}, see §\ref{sect:aut-vb-threezero}. We can also define automorphic line bundles
$\mathcal{L}'(\lambda)$ over $\FlKtr$.

\subsection{Weight-changing operators}
A great source of congruences are given by so-called \textit{weight-changing operators}.
These are Hecke-equivariant maps between spaces of automorphic forms in characteristic $p$.
There are two well-established types of weight-changing operators that have been
developed for Hodge-type Shimura varieties (which includes our setting).

First are the \textit{partial Hasse invariants}, which were defined
in \cite{flag-strata}. These are certain automorphic forms with trivial Hecke-action, so that multiplication by them provide weight-changing operators.
They are naturally defined
on the corresponding flag Shimura varieties.

Next are the \textit{theta operators}.
These are certain Hecke-equivariant differential operators between spaces
of automorphic forms of different weights.
They have been defined for Hodge-type Shimura varieties in \cite{ortiz}.

\subsection{Mod $p$ local-global compatibility}
Just like in the modular curve case, we can associate mod $p$ Galois representations to
mod $p$ Picard modular forms. More precisely, given a mod $p$ Picard modular form
$f$, we obtain a Galois representation $\rho_f : G_E \rightarrow \GL_3(\Fpb)$ that is essentially conjugate
self-dual.
We can now state \textit{mod $p$ local-global compatibility for irregular weight
Picard modular forms}.
\begin{conjecture}
  For any $f \in H^0(\Sh_{K, \Fpb}, \omega(k, 2))$ with $0 \le k \le p - 1$, the
  associated local Galois representation $\rho_f|_{G_{E_p}}$ admits an essentially conjugate self-dual crystalline lift
  with $\tau_0$-labeled Hodge--Tate cocharacter $\mu(\lambda) := (k + 2, 1, 1)$.\end{conjecture}
\begin{remark}
  \begin{enumerate}
  \item This can be seen as an analogue of Theorem
    \ref{thm:mod-p-LGC}, as (in the notation of that
    theorem) $\rho_f|_{G_{\Qp}}$ being unramified is equivalent to it admitting a crystalline lift
    of Hodge--Tate weights $\{0, 0\}$.
  \item We do not expect the analogue of the conjecture to hold for $k \ge p$.
    The reason is that for those weights, there
    exists a partial Hasse invariant giving an injective weight-changing operator from some regular weight
    to that irregular weight. But on the Emerton--Gee stack, the corresponding irregular locus will have codimension $1$, whereas
    the corresponding regular locus will have codimension $0$. In other words, we expect the space of
    irregular weight mod $p$ Picard modular forms to be too big in those cases.
  \item The explicit formula we use for the Hodge--Tate cocharacter draws from the characteristic $0$
    setting, where it is predicted that the infinitesimal character of an automorphic representation
    determines the Hodge--Tate weights of the associated Galois representation \cite[Remark
    3.2.3]{BG}.
  \end{enumerate}
\end{remark}

\subsection{Main result}
The main result of this paper is the following, establishing an exotic congruence
from irregular weight Picard modular forms on $\ShKtwo$ to $p$-restricted regular weight Picard modular forms on $\ShKtr$, where the regular weight can be pinned down to one out of two options. The following is a simplified version of Theorem \ref{thm:main-thm}
which avoids introducing additional notation at this stage (in particular, we do not account for determinant or multiplier twists).
\begin{theorem} \label{thm:introduction}
  Fix an irregular weight $\lambda = (k, 2)$ with $0 \le k \le p - 1$.
  Set $\lambda_1 = (p - 1 - k, k)$ and $\lambda_2 = (1, p - 2 - k)$.
  Then, if $k \le p - 2$ and $p > 3$, every mod $p$ Hecke eigensystem appearing in $H^0(\Sh_{K, \Fpb}, \omega(\lambda))$, must appear in
  either $H^0(\Sh'_{K', \Fpb}, \omega'(\lambda_1))$ or
  $H^0(\Sh'_{K', \Fpb}, \omega'(\lambda_2))$.
  If $k = p - 1$, then the corresponding mod $p$ Hecke eigensystem must appear
  in $H^0(\Sh'_{K', \Fpb}, \omega'(\lambda_1))$.
\end{theorem}
Figure \ref{fig:thm-congruences} illustrates the corresponding Hodge--Tate cocharacters $\mu(\lambda)$ and $\mu(\lambda_i)$, viewed
inside the two $p$-restricted alcoves. The yellow points indicates the corresponding Hodge--Tate weights
for $k = 0$, and the lines are directed, meaning that as $k$ increases, the Hodge--Tate weight moves
along the corresponding line.
Note that the red line, corresponding to $\mu(\lambda_2)$, is
of length $p - 1$ whereas the other two are of length $p$. This is because in the case $k = p - 1$, the theorem provides a single possible regular weight $\mu(\lambda_1)$.

Let us make some remarks about this result.
\usetikzlibrary{decorations.markings}
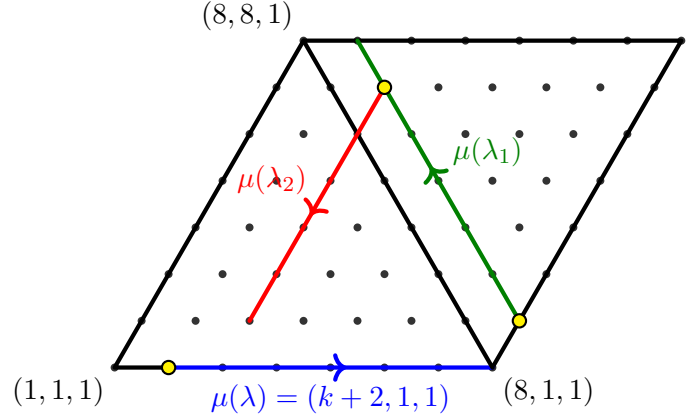
\begin{figure}
\centering
\begin{tikzpicture}[
    x={(0.714cm, 0cm)}, 
    y={(60:0.714cm)},
    mid arrow/.style={postaction={decorate, decoration={markings, mark=at position 0.55 with {\arrow{>}}}}}
]

    \foreach \a in {0,...,7} {
        \foreach \b in {0,...,7} {
            \fill[black!80] (\a,\b) circle (1.5pt);
        }
    }

    \draw[ultra thick, line join=round, line cap=round] (6,7) -- (0,7) -- (0,0) -- (1,0);
    \draw[ultra thick, line join=round, line cap=round] (7,0) -- (7,7) -- (6,7);
    
    \draw[ultra thick, blue, line cap=round, mid arrow] (1,0) -- (7,0) node[midway, below=2pt] {$\mu(\lambda) = (k + 2, 1, 1)$};
    
    \draw[ultra thick, line cap=round] (7,0) -- (0,7);

    \draw[ultra thick, green!50!black, line cap=round, mid arrow] (7,1) -- (1,7) node[midway, above right=1pt] {$\mu(\lambda_1)$};

    \draw[ultra thick, red, line cap=round, mid arrow] (2,6) -- (2,1);

    \path (2,6) -- (2,1) node[midway, above left, text=red] {$\mu(\lambda_2)$};

    \node[below left] at (0,0) {$(1, 1, 1)$};
    \node[below right] at (7,0) {$(8,1,1)$};
    \node[above left] at (0,7) {$(8,8,1)$};
    
    \fill[yellow, draw=black, thick] (1,0) circle (2.5pt); 
    \fill[yellow, draw=black, thick] (7,1) circle (2.5pt); 
    \fill[yellow, draw=black, thick] (2,6) circle (2.5pt); 

\end{tikzpicture}
\caption{Illustration of the Hodge--Tate weights in Theorem \ref{thm:introduction} for
$p = 7$ (up to determinant twists).}
\label{fig:thm-congruences}
\end{figure}

\begin{remark}
  \begin{enumerate}
    \item The case $k = p - 1$ is simultaneously the most important case of this theorem, and also when
      it is the most powerful. This is because, in this case, one can check that the partial Hasse invariants do not produce $p$-restricted congruences. We therefore expect that, in order to prove mod $p$ local-global compatibility in this case, one would need the congruence provided by our theorem, as well as a congruence produced by a theta operator.
    \item The case $k = 0$ is also important, but unfortunately provides no additional information. This is because
      $\lambda_2$ can already be obtained from a partial Hasse invariant. In upcoming work, we show
      how to eliminate $\lambda_2$ under genericity conditions on the Hecke eigensystem.
    \item For the cases $1 < k < p - 1$, the weight $\mu(\lambda_2)$ is Fontaine--Laffaille (i.e. lies
      in the lowest alcove). We expect that in those cases, one could eliminate this weight by a
      purely Galois-theoretic argument on the
      Emerton--Gee stack (using the partial Hasse and theta congruences available in
      these cases).
    \item Notice that, unlike the partial Hasse invariants and theta operators, the congruence we obtain is not given by a weight-shift but also involves reflections.
      A similar phenomenon happens in the setting of \cite{KLS}.
  \end{enumerate}
\end{remark}

\subsection{Geometric Jacquet--Langlands}

The geometric Jacquet--Langlands phenomenon in our setting was first worked out
by Vollaard \cite{vollaard} (although in a slightly different language). Namely, the paper completely describes the $1$-dimensional supersingular locus of
$\Sh_K$. In particular, it is reducible, and the irreducible components are all given by Fermat curves (i.e. curves in $\mathbb{P}^2$ cut out by the equation $x^{p + 1} + y^{p + 1} + z^{p + 1} = 0$). Furthermore, the set of irreducible components is in bijection with the Shimura set $\Sh_{K'}'$.
This can be rephrased as the existence of a correspondence
\[\begin{tikzcd}
	& {\mathcal{C}} & \\
	{\Sh_K} && {\Sh_{K'}',}
	\arrow["p"', from=1-2, to=2-1]
	\arrow["p'", from=1-2, to=2-3]
\end{tikzcd}\]
where $\mathcal{C}$ is given by a disjoint union of Fermat curves, one for each point in
$\Sh_{K'}'$. The map to $\Sh_K$ is then given by mapping each Fermat curve to the corresponding
irreducible component of the supersingular locus.
Let us note here that the notation $\mathcal{C}$ will only be used in the introduction (later there will be more precise notation relating it to a certain stratification).

This correspondence can be interpreted in a more sophisticated way on the level of flag Shimura varieties.
There, we can consider the \textit{coarse stratification}, which
generalizes the $a$-number stratification on the Shimura variety.
We can then view the above correspondence as the identification of $\mathcal{C}$ with certain strata
in the respective flag Shimura varieties $\Fl_K$ and $\FlKtr$. This gives us the following commutative diagram
\[\begin{tikzcd}
	&& {\mathcal{C}} && \\
	{\Sh_K} & {\Fl_K} && {\Fl'_{K'}} & {\Sh'_{K'},}
	\arrow["p"', bend right=14, from=1-3, to=2-1]
	\arrow["i"', hook', from=1-3, to=2-2]
	\arrow["{i'}", hook, from=1-3, to=2-4]
	\arrow["{p'}", bend left=14, from=1-3, to=2-5]
	\arrow["\pi"', from=2-2, to=2-1]
	\arrow["{\pi'}", from=2-4, to=2-5]
\end{tikzcd}\]
where $i$ and $i'$ are closed immersions.
\subsection{The exotic congruence}
The exotic congruence will be obtained as follows.
First, given an irregular weight $\lambda$, we
compute $i^*\mathcal{L}(\lambda)$ in terms of line bundles
on the flag Shimura variety $\FlKtr$.
More precisely, we find a weight $\lambda'$ such that
$i^*\mathcal{L}(\lambda) \cong i'^*\mathcal{L}'(\lambda')$.
Next, we must compute $p'_*i'^*\mathcal{L}'(\lambda')$ in terms of automorphic
vector bundles on $\ShKtr$, whose weights we can describe explicitly.
Pulling back a Picard modular form from $\FlKtwo$ to $\mathcal{C}$ and then pushing forward to $\ShKtr$
will then give us the desired exotic congruence.

Let us end the introduction with the following remark, regarding the possibility of eliminating the weight $\lambda_2$ without appealing to Galois theory.
An initial hope of ours was that, by working on the supersingular locus in $\ShKtwo$ instead of on $\mathcal{C}$ (which can be identified with its normalization), one should be able to eliminate $\lambda_2$ by analyzing the local intersections between the Fermat curves.
However, it seems that the obstruction rather comes from the ``global structure'' of the
supersingular locus, i.e. how the tree of Fermat curves from \cite{vollaard} folds in on itself. More precisely,
one can show that $\lambda_2$ can not be eliminated on the level of reduced Rapoport--Zink spaces.
In other words, on the reduced Rapoport--Zink space there will exist sections which obtain weight
$\lambda_2$ through an analogue of the exotic congruence on the level of reduced Rapoport--Zink spaces.
However, we do not expect such sections to descend to the supersingular locus of the Shimura variety.
Since we do not know how to get explicit control of this quotient, using it to eliminate $\lambda_2$ seems very difficult.

\subsection{Contents of the paper}
In §\ref{sect:dieudonne-main}, we start by recalling some Dieudonné theory from \cite{bbm}.
In particular, we are interested in how the Hodge filtration changes upon a quotient
by an $F, V$-kernel. This is described in Lemma \ref{lem:hodge-seq-final},
and will later be useful when computing the automorphic line bundles on $\FlKtwo$ restricted to $\mathcal{C}$.
Next, in §\ref{sect:shimura-main}, we introduce the
Shimura varieties $\ShKtwo$ and $\ShKtr$, their associated flag varieties and the
automorphic vector bundles. We also compute explicitly some partial Hasse invariants
in both cases, and write down explicit conditions defining certain strata in their
respective coarse stratifications.
In §\ref{sect:exotic-hecke-main} we establish a
version of the geometric Jacquet--Langlands phenomenon first discovered in \cite{vollaard}.
The main difference to Vollaard's work is that,
since we are interested in describing the automorphic line bundles,
we must also describe the universal abelian variety together with its Hodge filtration on this stratum.
The main new result of the section is Theorem \ref{thm:aut-vb-descr}, which gives a description of the automorphic line bundle
$\mathcal{L}(\lambda)$ in terms of automorphic line bundles on a partial flag Shimura variety.
In §\ref{sect:global-sections}, our goal is to compute the
pushforward of $\mathcal{L}(\lambda)|_{\mathcal{C}}$ to $\ShKtr$ in terms of automorphic vector
bundles with explicit weights. This is Theorem \ref{thm:main-decomp}.
We show how, using the stack of
$G$-zips, this reduces to the calculation of the cohomology of a certain line bundle on a Fermat curve, as a $G(\Fp)$-representation.
This calculation is then the content of Theorem \ref{thm:rep-decomposition}.
Finally, in §\ref{sect:applications} we apply Theorem \ref{thm:main-decomp} to obtain our main result, Theorem \ref{thm:main-thm}.
This requires us to show that restriction to the stratum $\mathcal{C}$ is injective on
global sections, which we do by an argument using partial Hasse invariants.

\subsection{Acknowledgements}
I want to thank my advisor George Boxer for suggesting this problem,
and for the numerous discussions we have had whilst I was working on this topic.
I also want to thank my other advisor, Fred Diamond, for comments and helpful discussions.
Finally, I want to also thank Martin Ortiz and Yicheng Yang for conversations about this work.

\subsection{Statement of AI use}
Gemini was used to produce Figure
\ref{fig:thm-congruences},
and it was also helpful when producing the notation section.
Other than this, no AI has been used in the creation of this article.

\subsection{Notation}
\subsubsection{Frobenius}
Given a scheme $S / \Fp$, let $\sigma_S$ (or simply $\sigma$) denote the absolute $p$-power Frobenius endomorphism on $S$. Given a scheme $X / S$, we set $X^{(p)} := X \times_{S, \sigma} S$. For a sheaf $\mathcal{F}$ over $S$, we will use the notation $\mathcal{F}^{(p^r)} := \sigma^{[r], *}_S \mathcal{F}$ to denote the $r$-fold pullback along $\sigma_S$.

\subsubsection{Witt vectors}
Given a perfect field $k/\Fp$, we let $W(k)$ denote the ring of Witt vectors of $k$. By an abuse of notation, we will also write $\sigma$ for the canonical lift of the Frobenius automorphism $\sigma_k$ to $W(k)$. Furthermore, we will use the notation $\Z_{p^2} := W(\Fpt)$ and $\Q_{p^2} := \Z_{p^2}[\frac{1}{p}]$.

\subsubsection{Characters and cocharacters}
Given an algebraic group $G$ over a field $k$,
we let $X^*(G)$ (resp.\ $X_*(G)$) denote the group of characters (resp.\ cocharacters) of $G_{\overline{k}}$.

\subsubsection{Root datum} Let $G$ be a connected reductive group equipped with a maximal torus $T \subseteq G$. Then we will let $R$ denote the associated set of $T$-roots in $X^*(T)$. We denote the set of simple roots by $\Delta$. Given a choice of a Borel subgroup containing $T$, we denote the corresponding set of positive roots by $R^+$ (and the positive simple roots by $\Delta^+$). We say a character $\lambda \in X^*(T)$ is $p$-restricted if $0 \le \langle \lambda, \alpha^\vee \rangle \le p - 1$ for all $\alpha \in \Delta^+$.

\subsubsection{Duals} We will employ various notations for duals. For a vector space $V$, $V^*$ will denote its dual vector space. For an abelian variety or a $p$-divisible group $A$, we let $A^\vee$ denote its dual. Finally, for a finite flat commutative group scheme $G$ over a base $S$, $G^D$ will denote its Cartier dual.

\subsubsection{Miscellaneous}
Given a morphism of schemes $X \xrightarrow{\pi} S$, we let $\mathcal{H}^i_{dR}(X/S)$ (or just $\mathcal{H}^i_{dR}(X)$) denote the
de Rham cohomology sheaves over $S$.
If $X$ is a group scheme, we let $\omega_{X/S}$ (often just $\omega_X$) denote the co-Lie algebra.
Similarly, we denote the Lie algebra by $\Lie(X/S)$ (often just $\Lie(X)$).

\section{Dieudonné theory} \label{sect:dieudonne-main}
In this section we will recall some results of \cite{bbm}, adapted for our purposes, and derive some
basic consequences. Throughout this section, we will fix a finite field $k/\Fp$ and a scheme
$S$ over $\Spec k$.

\subsection{Some remarks on the Crystalline site}

In \cite{bbm} they associate to any abelian sheaf $G$ on $S$, a sheaf
that we denote by $\mathcal{D}_{cris}(G)$ on the
crystalline site of $S$ over $W(k)$. This construction
is contravariantly functorial in $G$. Let us remark that whilst in \cite{bbm}, they use the notation
$\mathbb{D}(G)$ for this sheaf, we choose to use $\mathcal{D}_{cris}$ instead.
There is a morphism of
topoi $Sh((S/W(k))_{cris}) \rightarrow Sh(S_{zar})$ (Here $S_{zar}$ denotes the small Zariski site
on $S$), allowing us to pull back sheaves on the crystalline site to the Zariski site on $S$.
We can now produce the following two objects of interest from $\mathcal{D}_{cris}(G)$.
\begin{enumerate}
  \item Pulling back $\mathcal{D}_{cris}(G)$ to $S_{zar}$ gives us an
    $\mathcal{O}_S$-module $\mathcal{D}(G)_S$. We will call this the
    \textit{Dieudonné space} associated to $G$.
  \item When $S = \Spec k$, we let $\mathcal{D}(G)$ denote the $W(k)$-module obtained by taking
    global sections of $\mathcal{D}_{cris}(G)$ over the crystalline site. This is the
    Dieudonné module associated to $G$.
\end{enumerate}
The sheaf $\mathcal{D}_{cris}(G)$ always comes equipped with a Frobenius morphism
\[
  F : \mathcal{D}_{cris}(G)^{(p)} \longrightarrow \mathcal{D}_{cris}(G),
\]
This in turn induces Frobenius morphisms on $\mathcal{D}(G)_S$ resp. $\mathcal{D}(G)$ (for $S =
\Spec k$).
Similarly, whenever $G$ possesses a Verschiebung morphism, we obtain one on $\mathcal{D}_{cris}(G)$.
It will be denoted by $V$, and it will satisfy the standard compatibility with $F$. If any confusion
can occur, we will add a subscript to denote the associated group (i.e. $F_G$ resp. $V_G$).
\subsection{Dieudonné spaces}
The following proposition sums up the main properties that we will use about the
Dieudonné space associated to a $p$-divisible group.
\begin{proposition} \label{prop:dieu-space-props}
  Let $X/S$ be a $p$-divisible group. Then we have
  \begin{enumerate}[(i)]
    \item $\mathcal{D}(X)_S$ is locally free of rank equal to the height of $X$.
    \item $\mathcal{D}(X)_S \cong \mathcal{D}(X)_S / p\mathcal{D}(X)_S \cong \mathcal{D}(X[p])_S$, where $X[p]$ denotes the $p$-torsion subgroup.
    \item We have a canonical \textit{Hodge filtration}
    \[
      0 \longrightarrow \omega_{X} \longrightarrow \mathcal{D}(X)_S \longrightarrow \Lie(X^\vee) \longrightarrow 0.
    \]
  \item If there exists an abelian scheme $A/S$ such that $X \cong A[p^\infty]$, then we have
    \[
      \mathcal{D}(X)_S \cong \mathcal{D}(A)_S \cong \Hdr(A/S).
    \]
    Moreover, this isomorphism is compatible with Hodge filtrations.
    \item If $S = \Spec k$, we have $\mathcal{D}(X)_{\Spec k} \cong \mathcal{D}(X) / p \mathcal{D}(X)$.
  \end{enumerate}
\end{proposition}
\begin{proof}
  The first three parts follow from Théorème 3.3.10, Proposition 3.3.13 and Corollarie 3.3.5 of \cite{bbm} respectively.
  Next, iv) follows from Propositions 3.3.7 and 2.5.8.
  The final result follows from the fact that $\mathcal{D}(X)$ is a crystal.
\end{proof}
From now on, due to part (iv) of the above proposition, we will generally replace $\mathcal{D}(A)_S$ with
$\Hdr(A)$ for abelian schemes $A$ over $S$.
\begin{remark}
  Part (v) of the proposition shows that for $S = \Spec k$, the notion of ``Dieudonné space'' we use here coincides with
the one in \cite[(1.3)]{vollaard-wedhorn}.
\end{remark}
The following lemma will in particular allow us to describe $\mathcal{D}(G)_S$ for $G$ a finite flat group scheme
annihilated by $F_G$.
\begin{proposition} \label{prop:frob-quot}
  Given a finite flat group scheme $G$ over $S$ (or a $p$-divisible group),
  we have a canonical isomorphism
  \[
    \mathcal{D}(G)_S / F(\mathcal{D}(G)^{(p)}_S) \xlongrightarrow{\sim} \omega_G^{(p)}.
  \]
  Moreover, this isomorphism fits into the following commutative diagram,
  where the arrow on the right
  is obtained from the $(p)$-twisted Hodge filtration.
  \[\begin{tikzcd}
    {\mathcal{D}(G)_S / F(\mathcal{D}(G)^{(p)}_S)} & {\omega_G^{(p)}} \\
    & {\mathcal{D}(G)_S^{(p)}}
    \arrow["\sim", from=1-1, to=1-2]
    \arrow["V", from=1-1, to=2-2]
    \arrow[from=1-2, to=2-2]
  \end{tikzcd}
  \]
\end{proposition}
\begin{proof}
  This follows from \cite[Proposition 4.3.10]{bbm}, after pullback to $S_{zar}$.
\end{proof}

\subsection{Dieudonné complexes}
Given a finite flat group scheme $G$ over $S$, let $\Delta(G)_S \in D(\mathcal{O}_S)$ denote the (contravariant) \textit{Dieudonné complex}, which we define to be the pullback to the Zariski site of the complex defined in \cite[Définition 3.1.5]{bbm}.
Since inverse image functors between topoi are exact, it follows (from the definition)
that $\mathcal{H}^1(\Delta(G)_S) = \mathcal{D}(G)_S$.
The reason we are interested in this complex is that, as we shall soon see, certain theorems about Dieudonné modules of
$p$-divisible groups have analogue statements for Dieudonné complexes of finite flat group schemes.
The following theorem shows us that $\Delta(G)_S$ is not so complicated.
\begin{theorem}
  For a finite flat group scheme $G$, $\Delta(G)_S$ is representable by a perfect complex
  of $\mathcal{O}_S$-modules, with perfect amplitude contained in $[0, 1]$.
\end{theorem}
\begin{proof}
   This is \cite[Theorem 3.1.2 (i)]{bbm}.
\end{proof}
For a finite flat group scheme $S$, we will let $\ell_{G/S} \in \mathcal{D}(\mathcal{O}_S)$ denote its co-Lie complex.
Since finite flat group schemes are locally complete intersections,
$\ell_{G/S}$ has cohomological amplitude contained in the range $[-1, 0]$.
Recall that $\omega_G \cong \mathcal{H}^0(\ell_{G/S})$ and, following
\cite{bbm}, we will use the notation $n_G := \mathcal{H}^{-1}(\ell_{G/S})$ and $\nu_G :=
\mathcal{H}^1(\ell_{G/S}^*)$.

The Dieudonné complex admits a well-behaved ``Hodge
filtration'' \cite[Theorem 3.2.10]{bbm}. We recall the statement in our notation here.
\begin{theorem} \label{thm:ff-hodge}
  Fix a finite flat group scheme $G$ over $S$, and let $G^D/S$ denote its Cartier dual. There exists the following exact triangle
  in $D(\mathcal{O}_S)$ (functorial in $G$)
  \[
    \ell_{G/S}[-1] \longrightarrow \Delta(G)_S \longrightarrow \ell_{G^D/S}^{*} \xrightarrow{+1}.
  \]
  In particular, we have the following long exact sequence
  \[
    0 \longrightarrow n_G \longrightarrow \mathcal{H}^0(\Delta(G)_S) \longrightarrow \Lie(G^D) \longrightarrow
    \omega_G \longrightarrow \mathcal{D}(G)_S \longrightarrow \nu_{G^D} \longrightarrow 0.
  \]
\end{theorem}
The Dieudonné complex and its Hodge filtration behaves well with respect to isogenies.
\begin{lemma} \label{lem:hodge-seq}
  Given an isogeny of $p$-divisible groups $f : X_0 \rightarrow X_1$ over $S$, let $G$ denote the
  kernel. Then we have the following isomorphism of perfect complexes
  \[
    \Delta(G)_S \cong [\mathcal{D}(X_1)_S \xlongrightarrow{-\mathcal{D}(f)} \mathcal{D}(X_0)_S].
  \]
  This isomorphism is also compatible with Hodge filtrations, and in particular
  we get the following commutative diagram, where all rows and columns are
  exact
  \[\begin{tikzcd}
    && 0 & 0 & 0 && \\
    & 0 & {n_G} & {\omega_{X_1}} & {\omega_{X_0}} & {\omega_G} & 0 \\
    & 0 & {\mathcal{H}^0(\Delta(G)_S)} & {\mathcal{D}(X_1)_S} & {\mathcal{D}(X_0)_S} & {\mathcal{D}(G)_S} & 0 \\
    {} & 0 & {\Lie (G^D)} & {\Lie(X_1^\vee)} & {\Lie(X_0^\vee)} & {\nu_{G^D}} & 0 \\
    &&& 0 & 0 & 0
    \arrow[from=1-3, to=2-3]
    \arrow[from=1-4, to=2-4]
    \arrow[from=1-5, to=2-5]
    \arrow[from=2-2, to=2-3]
    \arrow[from=2-3, to=2-4]
    \arrow[from=2-3, to=3-3]
    \arrow[from=2-4, to=2-5]
    \arrow[from=2-4, to=3-4]
    \arrow[from=2-5, to=2-6]
    \arrow[from=2-5, to=3-5]
    \arrow[from=2-6, to=2-7]
    \arrow[from=2-6, to=3-6]
    \arrow[from=3-2, to=3-3]
    \arrow[from=3-3, to=3-4]
    \arrow[from=3-3, to=4-3]
    \arrow[from=3-4, to=3-5]
    \arrow[from=3-4, to=4-4]
    \arrow[from=3-5, to=3-6]
    \arrow[from=3-5, to=4-5]
    \arrow[from=3-6, to=3-7]
    \arrow[from=3-6, to=4-6]
    \arrow[from=4-2, to=4-3]
    \arrow[from=4-3, to=4-4]
    \arrow[from=4-4, to=4-5]
    \arrow[from=4-4, to=5-4]
    \arrow[from=4-5, to=4-6]
    \arrow[from=4-5, to=5-5]
    \arrow[from=4-6, to=4-7]
    \arrow[from=4-6, to=5-6]
  \end{tikzcd}\]
\end{lemma}
\begin{proof}
  This follows from \cite[Proposition 3.3.13]{bbm}.
\end{proof}
\subsection{Quotients by $F, V$-kernels} \label{sect:ht-one}
We say that a finite flat group scheme $G/S$ is an $F, V$-kernel, if it
is annihilated by both $F_G$ and $V_G$.
In this section we aim to understand how the Dieudonné spaces of $p$-divisible
groups (and their Hodge filtrations) change upon quotienting by an $F, V$-kernel.
\begin{lemma} \label{lem:hodge-seq-FV}
  Given an $F, V$-kernel $G$ over $S$, the long exact sequence of Theorem \ref{thm:ff-hodge} becomes
  \[
    0 \longrightarrow n_G \xlongrightarrow{\sim} \mathcal{H}^0(\Delta(G)_S) \xlongrightarrow{0} \Lie(G^D) \xlongrightarrow{\sim}
    \omega_G \xlongrightarrow{0} \mathcal{D}(G)_S \xlongrightarrow{\sim} \nu_{G^D} \longrightarrow 0.
  \]
\end{lemma}
\begin{proof}
    This follows from \cite[Lemme 4.3.2]{bbm}, which implies that all terms
    of the long exact sequence are locally free of the same rank.
\end{proof}
For the statement of the next lemma, the following notation will be convenient.
Given a morphism $F : \mathcal{E} \to \mathcal{F}$ of vector bundles over some scheme $S$,
and a subquotient $\mathcal{Q}$ of $\mathcal{E}$ presented by
\[
  0 \longrightarrow \mathcal{E}'' \longrightarrow \mathcal{E}' \longrightarrow \mathcal{Q} \longrightarrow 0,
\]
where $\mathcal{E}'' \subseteq \mathcal{E}' \subseteq \mathcal{E}$ are subsheaves.
Then we will let $F(\mathcal{Q})$ denote $F(\mathcal{E}') / F(\mathcal{E}'')$, which is 
a (possibly trivial) subquotient of $\mathcal{F}$.

The following lemma will allow us to describe the action of Frobenius and
Verschiebung on the quotient of a $p$-divisible group by an $F, V$-kernel.

\begin{lemma} \label{lem:frob-ver-descr}
  Let $S$ be a scheme of characteristic $p$, and let $\mathcal{X} \rightarrow \mathcal{Y}$
  be an isogeny of $p$-divisible groups over $S$, with kernel $G$.
  Consider the exact sequence
  \[
    0 \longrightarrow \mathcal{H}^0(\Delta(G)_S) \longrightarrow \mathcal{D}(\mathcal{Y})_S \xlongrightarrow{\Phi} \mathcal{D}(\mathcal{X})_S \longrightarrow \mathcal{D}(G)_S \longrightarrow 0.
  \]
  If $G$ is an $F, V$-kernel, then both $V_{\mathcal{X}}(\mathcal{D}(G)_S)$
  and $F_{\mathcal{X}}(\mathcal{D}(G)_S^{(p)})$ are non-zero, lying in
  the image of $\Phi$, and we have that
  \[
    F_{\mathcal{Y}}(V_{\mathcal{X}}(\mathcal{D}(G)_S)) = \mathcal{H}^0(\Delta(G)_S)
  \]
  and
  \[
    V_{\mathcal{Y}}(F_{\mathcal{X}}(\mathcal{D}(G)_S^{(p)})) = \mathcal{H}^0(\Delta(G)_S)^{(p)}.
  \]
  In particular, $\mathcal{D}(G)_S \cong \mathcal{H}^0(\Delta(G)_S)$.
\end{lemma}
\begin{proof}
  We only prove the first case, the second is identical.
  First, the fact that $F_G$ is trivial implies that $\mathcal{D}(\mathcal{X})_S \rightarrow \mathcal{D}(G)_S$
  factors through $\mathcal{D}(\mathcal{X})_S / \Img F_{\mathcal{X}} \rightarrow \mathcal{D}(G)_S$.
  Next, we use that $\Img F_{\mathcal{X}} = \ker V_{\mathcal{X}}$ in $\mathcal{D}(\mathcal{X})_S$ (since $\mathcal{X}[p]$ is a $BT_1$)
  to deduce that $\mathcal{D}(G)_S \xrightarrow{V_{\mathcal{X}}} V_{\mathcal{X}}(\mathcal{D}(G)_S)$ is an isomorphism.
  Hence,
  $\mathcal{E} := V_{\mathcal{X}}(\mathcal{D}(G)_S)$ is a vector bundle of rank equal to that of $\mathcal{D}(G)_S$, and in particular also to that of $\mathcal{H}^0(\Delta(G)_S)$.
  Since $V_G$ is trivial, we know that $\mathcal{E}$ must lie in the image of $\Phi$.
  Hence, we can view $\mathcal{E}$ as a subquotient of $\mathcal{D}(\mathcal{Y})_S$ to make sense of the expression
  $F_{\mathcal{Y}}(\mathcal{E})$.
  It remains to see that this is equal to $\mathcal{H}^0(\Delta(G)_S)$.
  First, let us show that $F_{\mathcal{Y}}(\mathcal{E})$ has the same rank as $\mathcal{E}$.
  To do this, we need to see that $\mathcal{E}$ does not lie in the image of $V_{\mathcal{Y}}$.
  But this is impossible, since $\mathcal{D}(G)_S$ does not lie
  in the image of $\Phi$.
  Thus, $F_{\mathcal{Y}}(\mathcal{E})$ has rank equal to that of $\mathcal{E}$, but
  since $F_{\mathcal{X}}(\mathcal{E}) = 0$, we furthermore have that $F_{\mathcal{Y}}(\mathcal{E})
  \subseteq \ker \Phi = \mathcal{H}^0(\Delta(G)_S)$.
  Since this is a containment of subbundles of the same rank, it must be an
  equality.
\end{proof}
Now we can completely describe how the Hodge filtration
changes upon taking the quotient by an $F, V$-kernel. This lemma
will be repeatedly used in §\ref{sect:exotic-hecke-main}. 
\begin{lemma} \label{lem:hodge-seq-final}
  Let $S$ be a scheme of characteristic $p$, and let $\mathcal{X} \rightarrow \mathcal{Y}$
  be an isogeny of $p$-divisible groups over $S$, with kernel $G$.
  If $G$ is an $F, V$-kernel, then
  the associated short exact sequence of Hodge filtrations from Lemma \ref{lem:hodge-seq} can be described as
\[\begin{tikzcd}
	0 & {\omega_{G}^{(p)}} & {\omega_{\mathcal{Y}}} & {\omega_{\mathcal{X}}} & {\omega_{G}} & 0 \\
	0 & {\mathcal{H}^0(\Delta(G)_S)} & {\mathcal{D}(\mathcal{Y})_S} & {\mathcal{D}(\mathcal{X})_S} & {\mathcal{D}(G)_S} & 0 \\
	0 & {\omega_{G}} & {\Lie(\mathcal{Y}^\vee)} & {\Lie(\mathcal{X}^\vee)} & {\omega_{G}^{(p)}} & {0.}
	\arrow[from=1-1, to=1-2]
	\arrow[from=1-2, to=1-3]
	\arrow["{\simeq}", from=1-2, to=2-2]
	\arrow[from=1-3, to=1-4]
	\arrow[from=1-3, to=2-3]
	\arrow[from=1-4, to=1-5]
	\arrow[from=1-4, to=2-4]
	\arrow[from=1-5, to=1-6]
	\arrow["0", from=1-5, to=2-5]
	\arrow[from=2-1, to=2-2]
	\arrow[from=2-2, to=2-3]
	\arrow["0", from=2-2, to=3-2]
	\arrow[from=2-3, to=2-4]
	\arrow[from=2-3, to=3-3]
	\arrow[from=2-4, to=2-5]
	\arrow[from=2-4, to=3-4]
	\arrow[from=2-5, to=2-6]
	\arrow["{\simeq}", from=2-5, to=3-5]
	\arrow[from=3-1, to=3-2]
	\arrow[from=3-2, to=3-3]
	\arrow[from=3-3, to=3-4]
	\arrow[from=3-4, to=3-5]
	\arrow[from=3-5, to=3-6]
\end{tikzcd}\]
Here, the isomorphism $\mathcal{D}(G)_S \cong \omega_G^{(p)}$ is induced from $V_{\mathcal{X}} :
\mathcal{D}(\mathcal{X})_S \twoheadrightarrow \omega_{\mathcal{X}}^{(p)} \twoheadrightarrow \omega_{
G}^{(p)}$.
Furthermore, by viewing $\omega_G$ as a subquotient of $\mathcal{D}(\mathcal{Y})_S$ (using that the top
right vertical morphism is $0$), we get that $F_{ \mathcal{Y}}$ induces an isomorphism
\[
  F_{\mathcal{Y}} : \omega_{G}^{(p)} \xrightarrow{\sim} \mathcal{H}^0(\Delta(G)_S).
\]
\end{lemma}
\begin{proof}
  Applying Lemma \ref{lem:hodge-seq} and Lemma
\ref{lem:hodge-seq-FV} gives us the following diagram 
\[\begin{tikzcd}
  0 & {n_G} & {\omega_{\mathcal{Y}}} & {\omega_{\mathcal{X}}} & {\omega_{G}} & 0 \\
	0 & {\mathcal{H}^0(\Delta(G)_S)} & {\mathcal{D}(\mathcal{Y})_S} & {\mathcal{D}(\mathcal{X})_S} & {\mathcal{D}(G)_S} & 0 \\
	0 & {\omega_G} & {\Lie(\mathcal{Y}^\vee)} & {\Lie(\mathcal{X}^\vee)} & {\nu_{G^D}} & {0.}
	\arrow[from=1-1, to=1-2]
	\arrow[from=1-2, to=1-3]
	\arrow["\simeq", from=1-2, to=2-2]
	\arrow[from=1-3, to=1-4]
	\arrow[from=1-3, to=2-3]
	\arrow[from=1-4, to=1-5]
	\arrow[from=1-4, to=2-4]
	\arrow[from=1-5, to=1-6]
	\arrow["0", from=1-5, to=2-5]
	\arrow[from=2-1, to=2-2]
	\arrow[from=2-2, to=2-3]
	\arrow["0", from=2-2, to=3-2]
	\arrow[from=2-3, to=2-4]
	\arrow[from=2-3, to=3-3]
	\arrow[from=2-4, to=2-5]
	\arrow[from=2-4, to=3-4]
	\arrow[from=2-5, to=2-6]
	\arrow["\simeq", from=2-5, to=3-5]
	\arrow[from=3-1, to=3-2]
	\arrow[from=3-2, to=3-3]
	\arrow[from=3-3, to=3-4]
	\arrow[from=3-4, to=3-5]
	\arrow[from=3-5, to=3-6]
\end{tikzcd}\]
It remains to identify $\mathcal{D}(G)_S$ and $\mathcal{H}^0(\Delta(G)_S)$ (then we also identify
$n_G$ and $\nu_{G^D}$ for free).
By Proposition \ref{prop:frob-quot} , since $G$ is killed by $F_G$, we have that $\mathcal{D}(G)_S \cong \omega_{G}^{(p)}$.
Moreover, the second part of the proposition applied to the group scheme $\mathcal{X}[p]$ implies that this isomorphism is induced by
$V_{\mathcal{X}}$.
Now Lemma \ref{lem:frob-ver-descr} gives us that $F_{\mathcal{Y}}$ induces an
isomorphism $\omega_G^{(p)} \cong \mathcal{H}^0(\Delta(G)_S)$ (by viewing $\omega_G^{(p)}$ as a subquotient of $\mathcal{D}(\mathcal{Y})_S^{(p)}$).
In particular, we get that $n_G \cong \mathcal{H}^0(\Delta(G)_S) \cong \omega_G^{(p)}$.
This finishes the proof of the theorem.
\end{proof}

\section{Unitary Shimura varieties} \label{sect:shimura-main}
\subsection{Unitary groups} \label{sect:unitary-grp}
Fix an imaginary quadratic extension $E$ of $\Q$, such that $p$ is inert in $E$.
We will use $x \mapsto \overline{x}$ to denote the nontrivial automorphism of $E/\Q$.

Let $V$ be a $3$-dimensional $E$-vector space, equipped with a non-degenerate $E$-valued skew-Hermitian
pairing $\{ \cdot, \cdot \} : V \times V \rightarrow E$.
By taking traces, we obtain a non-degenerate $\Q$-valued alternating pairing $\langle \cdot,\cdot \rangle : V \times V \rightarrow \Q$ given by $\langle v, w \rangle := \Tr_{E/\Q} (\{v, w\})$.
Furthermore, this pairing satisfies $\langle xv,w \rangle = \langle v,\overline{x}w \rangle$ for all $x \in E$ and $v,w \in V$.
Using this, we can define an algebraic group $G$ over $\Q$,
which for a $\Q$-algebra $R$, has its group of $R$-points given by
\[
	G(R) = \{ (g, c) \in (\Res_{E/\Q} GL(V) \times \Gm)(R) \mid \langle gv,gw \rangle = c \langle v,w \rangle \}.
\]
Let us once and for all fix a basis $V \cong E^{\oplus 3}$. The skew-Hermitian space $V_\R \cong \C^3$ will have some signature
$(r, s)$, where $r + s = 3$, so in particular $G_{\R} \cong \GU(r, s)$.

As we will see, the reflex field of our Shimura variety will be given by $E$.
It is nonetheless good practice to separate the reflex field from the endomorphism structure.
Therefore, let $E_0$ denote another copy of $E$, which we will use to denote the reflex field (but we do not fix an isomorphism $E_0 \cong
E$). Fix a labeling $\Hom_{\Q}(E, E_0) = \Hom_{\Z}(\mathcal{O}_E, \mathcal{O}_{E_0}) = \{ \tau_0, \tau_1 \}$.
For any $\mathcal{O}_{E_0}$-algebra $R$, we will also let $\tau_i$ denote the corresponding map from $\mathcal{O}_E$.

We get a decomposition $V \otimes_{\Q} E_0 \cong V_{\tau_0} \oplus V_{\tau_1}$,
where $E$ acts on the right-hand side by $\tau_i$ on $V_{\tau_i}$.
Throughout this paper we will interpret the indices of $\tau_i$ as lying in $\Z/2\Z$.
Let $\lambda_V : V \xrightarrow{\sim} V^*$ denote the morphism associated to the pairing $\langle \cdot,
\cdot \rangle$.
Note that $\lambda_V \otimes_{\Q} E_0$ exchanges isotypic components, and hence we can consider its restriction to $V_{\tau_0}$
\[
  \lambda_{V, \tau_0} : V_{\tau_0} \xlongrightarrow{\sim} V_{\tau_1}^*.
\]
This recovers the skew-Hermitian pairing $\{ \cdot , \cdot \}$.
Finally, we have an identification
$G_{E_0} \cong \GL(V_{\tau_0}) \times \mathbb{G}_{m, E_0}$, since the action on $V_{\tau_1}$ can be determined from the action on $V_{\tau_0}$ via $\lambda_{V, \tau_0}$.

\subsubsection{Pairings in the type $(2, 1)$ and $(3, 0)$ cases}
We will mainly be interested in the two cases $(r, s) \in \left\{ (2, 1), (3, 0) \right\}$.
Let us for future use fix the pairing $\{ \cdot, \cdot \}$ in each case.
Fix an element $\zeta \in \mathcal{O}_E$ satisfying $\overline{\zeta} = -\zeta$.
In the $(2, 1)$-case, we will use the skew-Hermitian pairing given by
\[
  \{v,w\} := \overline{v}^t 
  \begin{pmatrix}
    0 & 0 & \zeta \\
    0 & \zeta & 0 \\
    \zeta & 0 & 0 \\
  \end{pmatrix}
  w.
\]
We make this choice so that the upper-triangular Borel descends from $G_{E_0}$ to $G$.
In the $(3, 0)$-case, we will use the standard skew-Hermitian pairing
\[
  \{v,w\}' := \overline{v}^t 
  \begin{pmatrix}
    \zeta & 0 & 0 \\
    0 & \zeta & 0 \\
    0 & 0 & \zeta \\
  \end{pmatrix}
  w.
\]
\subsubsection{An integral model}
Consider the lattice $\Lambda := \mathcal{O}_{E_p}^{\oplus 3} \subseteq V \otimes_{\Q} \Qp$ (with
respect to our fixed basis).
We will assume that the pairing $\{ \cdot, \cdot \}$ restricts to a perfect pairing on $\Lambda$ (which it does in the explicit cases considered above), and hence
$G$ admits a reductive model $\mathcal{G}$ over $\Zp$.
We will use the notation $G(\Zp) := \mathcal{G}(\Zp)$, and
$G_{\Fp}$ for the reductive special fiber of $\mathcal{G}$.
Let $W := \Lambda / p \Lambda$ denote the mod $p$ reduction, and note that
we can identify $G_{\Fp} = \GU(W)$ (where $W$ obtains an induced pairing from $\Lambda$).

\subsubsection{The maximal torus}
Let $T \subset G$ denote the maximal torus consisting of diagonal matrices with respect to our fixed
basis of $V$.
We can identify $T_{E_0} \cong \mathbb{G}_m^{\oplus 3} \times \mathbb{G}_m \subseteq
\GL(V_{\tau_0}) \times \mathbb{G}_m$, so it is in particular split over $E_0$ and we can identify $X^*(T) = \Z^{\oplus 3} \times \Z$.
Since $T$ is not split, $X^*(T)$ comes with a non-trivial Galois action (of the group $\Gal(E_0 /
\Q)$). One checks that the non-trivial element acts via
\[
  (a, b, c, d) \longmapsto (- c, - b, - a, a + b + c + d).
\]

\subsubsection{The Borel subgroup}
Using that $G_{E_0} \cong \GL(V_{\tau_0}) \times \mathbb{G}_{m, E_0}$, we let
$B \subseteq G_{E_0}$ denote the upper-triangular Borel subgroup.
Note that in the
type $(2, 1)$-case, our choice of pairing $\{ \cdot, \cdot\}$ implies that this
Borel descends to $G/\Spec \Q$.

\subsubsection{The root datum}
Let $R$ denote the set of roots of $G$ with respect to the torus $T$.
The choice of Borel $B$ induces a set of positive roots $R^+ \subseteq R$.
We have the two simple (positive) roots $\alpha_1 = (1, -1, 0, 0)$ and $\alpha_2 = (0, 1, -1, 0)$.
Let $\rho := \frac{1}{2}\sum_{\alpha \in R^+} \alpha$ denote the half-sum of positive roots.

\subsubsection{The Weyl group}
Let $W(T)$ be the Weyl group acting on $T$. It is isomorphic to the symmetric group
$S_3$, acting by permuting the basis elements of $T$. We will let $w_0$ denote the
longest element of the Weyl group, and we let $s_i$ denote the simple reflection corresponding to
$\alpha_i$. In terms of permutations, we can describe these as $w_0 = (13), s_1 = (12), s_2 = (23)$.
Later on, once we fix a Hodge cocharacter, we will have a Levi subgroup $M \subseteq G_{E_0}$ (containing $T$) which depends on the signature $(r, s)$.
We will let $w_c$ denote the longest element of the Weyl group of $M$, viewed inside $W(T)$.

\subsection{Unitary Shimura varieties}
Throughout this section we work with a fixed general type $(r, s)$ (with $r + s = 3$).

\subsubsection{The integral Shimura datum}
Let $\mathbb{S}$ denote the Deligne torus $\Res_{\C / \R} \Gm$. We can then
consider the homomorphism $h : \mathbb{S} \rightarrow G_\R$ given by $h(z) = \operatorname{diag}(z\mathbf{1}_{r},
\overline{z}\mathbf{1}_s)$.
Combining this with what was done in §\ref{sect:unitary-grp}, we now have an
integral PEL Shimura datum given by
$(E, V, \langle \cdot, \cdot \rangle, h, \Lambda, \mathcal{O}_{E_p})$.

\subsubsection{The Hodge cocharacter}
Associated to $h$, we can attach the Hodge cocharacter $\mu : \mathbb{G}_{m, E_0} \rightarrow G_{E_0}$ given by $\mu(z) = \operatorname{diag}(z \mathbf{1}_r, \mathbf{1}_s) \times z$. Its field of definition is $E_0$, which will be the reflex field of our Shimura varieties.
The stabilizer of $\mu$ will be a Levi subgroup which we denote by $M$.
From $\mu$ we also obtain a parabolic $P$, given by
\[
P = \{ g \in G_{E_0} \mid \lim_{t \to 0} \mu(t)g\mu(t)^{-1} \text{ exists} \}.
\]
$P$ will be isomorphic to the stabilizer of the corresponding Hodge filtration at any point of our Shimura variety.
We will let $P^-$ denote the opposite parabolic.
Both $P$ and $M$ can be defined integrally, and we let $P_{\Fpt}$ resp. $M_{\Fpt}$ denote their special fibers.
Note that the bulk of this paper will be taking place in characteristic $p$, so
it will be common to replace $P$ and $M$ with their special fibers to lighten notation (we will make it clear when this applies).
\subsubsection{The moduli space}
We now define the moduli space corresponding to our integral PEL datum.
Fix an isomorphism $\mathcal{O}_{E_{0, p}} \cong \Zpt$.
Our Shimura varieties will be defined integrally, as schemes over $\Spec \Zpt$, and with \textit{hyperspecial} level.
Fix an open compact subgroup $K^p \subseteq \adfp$, where $\adfp$ denotes the finite adeles
away from $p$. Let $K := K^pK_p$, where $K_p = G(\Zp)$ (i.e. the level is hyperspecial at $p$).

Let $\Sh_K$ denote the moduli problem sending a $\Zpt$-scheme $S$
to the set of equivalence classes $\{ (A, \iota, \lambda, \eta) \} / \sim$ where
\begin{enumerate}[(i)]
  \item $A \rightarrow S$ is an abelian scheme.
  \item The dual abelian scheme $A^\vee$ exists, and $\lambda : A \rightarrow A^\vee$ is a prime-to-$p$ quasi-polarization (meaning, it is
    invertible as an element of $\operatorname{Isog}(A, A^\vee) \otimes \Z_{(p)}$).
  \item $\iota : \mathcal{O}_{E} \rightarrow \End(A) \otimes \Z_{(p)}$
    is a homomorphism such that for all $x \in \mathcal{O}_{E}$, we have
   $\iota(\overline{x}) =\lambda^{-1} \circ \iota(x)^\vee \circ \lambda$.
   Furthermore, it needs to satisfy the signature condition $(r, s)$, i.e. for all $x \in \mathcal{O}_{E}$
    we have
    \[
      \operatorname{charpol}(\iota(x)|_{\Lie(A)}) = (T - \tau_0(x))^r(T - \tau_1(x))^s.
    \]
  \item $\eta$ is a $K^p$-level structure. By this we mean that, for each connected component $S_i$ of $S$,
    after fixing a geometric point $s_i \in S_i$, we have some fixed $K^p$-orbit
    \[
      \eta_i \in K^p \backslash \Isom_{G(\adfp)}(V \otimes_{\Q} \adfp, V^p(A_{s_i})).
    \]
    Here, $V^p(A_{s_i})$ denotes the rational Tate module of $A_{s_i}$ away from $p$.
    It comes equipped with the Weil pairing, taking values in $\adfp(1)$.
    The notation $\Isom_{G(\adfp)}$ means that we require our isomorphisms to
    respect the $\mathcal{O}_E$-action, and that there is some $\nu(\eta_i)$ making the following
    diagram commute, where the horizontal arrows are the pairings
  \[\begin{tikzcd}
    {V \otimes_{\Q} \adfp \times V \otimes_{\Q} \adfp} & \adfp \\
    {V^p(A_{s_i}) \times V^p(A_{s_i})} & {\adfp(1)}
    \arrow[from=1-1, to=1-2]
    \arrow["{\eta_i \times \eta_i}", from=1-1, to=2-1]
    \arrow["{\nu(\eta_i)}", from=1-2, to=2-2]
    \arrow[from=2-1, to=2-2]
  \end{tikzcd}\]
  Finally, we require that the $K^p$-orbit $\eta_i$ is $\pi_1(S_i, s_i)$-stable.
\end{enumerate}
By equivalence classes, we mean that we identify tuples $(A_i, \iota_i, \lambda_i, \eta_i)$ for $i =
0, 1$, if there exists a prime-to-$p$ quasi-isogeny $A_1 \dashrightarrow A_2$ preserving the
$\mathcal{O}_E$ action, and the $K^p$-level structure, as well as preserving
the polarization $\lambda$ up to a $\operatorname{Cont}(S, \Z_{(p)}^\times)$-scalar.

Note that (iii) implies that the relative dimension of $A$ must be $3$.

For small enough $K^p$, this moduli problem is representable by a smooth, quasi-projective scheme over $\Spec \Zpt$ of relative dimension $rs$ \cite[\S 5]{kottwitz}.
In the case $rs = 0$, we will refer to $\Sh_K$ as a \textit{Shimura set}, and
otherwise we will refer to it as a Shimura variety.
\subsubsection{Adelic double quotients}
Since $G$ satisfies the Hasse principle\footnote{By \cite[§7]{kottwitz} this reduces to $\Res_{E/\Q} \mathbb{G}_m$, for which
it holds by the Hasse norm theorem.},
the complex points of $\Sh_K$ admit a description as an adelic double quotient
\[
  \Sh_K(\C) = G(\Q) \backslash (X \times G(\adf) / K).
\]
Here $X$ denotes the $G(\R)$-conjugacy class of $h$. Notice that for $rs = 0$,
$X$ is just a point.
\subsection{Automorphic vector bundles} \label{sect:aut-vb}
Throughout this section, we will be working over the base $\Spec \Fpt$,
as our results only concern automorphic vector bundles in characteristic $p$ anyway.
To avoid cluttering the notation too much, we will therefore consider all
schemes to be implicitly base changed to $\Spec \Fpt$.
In particular, $T$ will be split in this section.
\subsubsection{Construction}
Let $(\mathcal{A}, \iota, \lambda, \eta)$ denote the universal family over $Sh_K$.
It admits a Hodge filtration, which is determined by its $\tau_0$-isotypic
component
\[
  0 \longrightarrow \omega_{\mathcal{A}, \tau_0} \longrightarrow 
  \Hdr(\mathcal{A})_{\tau_0} \longrightarrow
  \Lie(\mathcal{A}^\vee)_{\tau_0} \longrightarrow 0.
\]
Consider the space of trivializations
\[
  G_{dR} := \underline{\Isom}_{\Sh_K}(\Hdr(\mathcal{A})_{\tau_0},
  W_{\tau_0} \otimes_{\Fpt} \mathcal{O}_{\Sh_K}).
\]
By \cite[Theorem 3.4.1]{zhang}, $G_{dR} \rightarrow \Sh_K$ is a (left) étale $G$-torsor.
It admits a natural morphism to the flag space $G / P$, given by transferring the Hodge filtration
via the given trivialization.
This gives us a pullback square
\[\begin{tikzcd}
	{G_{dR}} & {G / P} \\
	{\Sh_K} & {\left[* / P\right].}
	\arrow[from=1-1, to=1-2]
	\arrow[from=1-1, to=2-1]
	\arrow[from=1-2, to=2-2]
	\arrow[from=2-1, to=2-2]
\end{tikzcd}\]
Let $\Rep(P)$ denote the category of finite dimensional algebraic $P$-representations, and for an algebraic stack $\mathfrak{X}$, we will let $VB(\mathfrak{X})$ denote the corresponding category of finite rank vector bundles. Then we get a functor from the above square via pullback
\[
  \omega : \Rep(P) \cong VB(\left[* / P \right]) \longrightarrow VB(\Sh_K).
\]
The image of this functor will be the \textit{automorphic vector bundles} on $\Sh_K$.

We are particularly interested in automorphic vector bundles induced from the
 torus $T \subseteq M$. Namely, given a character $\lambda \in X^*(T)$,
we get the induced representation $\Ind_{M \cap B}^{M} w_c \lambda$ of the Levi $M$,
which we can then inflate to a representation of $P$. This gives us a map which we
again denote by
\[
  \omega : X^*(T) \longrightarrow VB(\Sh_K).
\]
Given $\lambda \in X^*(T)$ we will say that $\omega(\lambda)$ is the automorphic
vector bundle of \textit{weight} $\lambda$.
The weights of the form $\lambda = (a, b, b - 2, d)$ with $a \ge b$ are the \textit{irregular} weights (see e.g. \cite[§2]{goldring} \cite[§3]{harris} \cite{BHR}).

Let us also introduce the notation $\delta := \omega(0, 0, 0, 1)$, and $\det := \omega(1, 1, 1, 0)$, corresponding to the multiplier resp. the determinant
characters of $G$ (restricted to $P$).
They are trivial as line bundles, but come with non-trivial Hecke actions.

\subsection{Unitary flag Shimura varieties}
We now introduce the \textit{flag Shimura varieties}. They are useful as the cohomology of automorphic vector bundles on $\Sh_K$ can be identified with the cohomology of certain automorphic line bundles on the flag Shimura variety.
For simplicity, we will continue working over $\Spec \Fpt$ in this section.

Fix a parabolic $Q \subseteq P$, then we define $\Fl_{K, Q}$ as the pullback
\[\begin{tikzcd}
	{\Fl_{K, Q}} & {\left[* / Q\right]} \\
	{\Sh_K} & {\left[* / P\right].}
	\arrow[from=1-1, to=1-2]
	\arrow["\pi_Q", from=1-1, to=2-1]
	\arrow[from=1-2, to=2-2]
	\arrow[from=2-1, to=2-2]
\end{tikzcd}\]
If $Q$ is not a Borel, we will refer to these as \textit{partial flag Shimura varieties}. From now on set $Q = B$ (we will return to $Q \neq B$ in §\ref{sect:partial-flag}).
In this case, we will omit the $B$-subscript from the notation.
The variety $\Fl_K$ admits a moduli interpretation as follows.
Given a scheme $S$ over $\Spec \Fpt$, $\Fl_K(S)$
can be described as the set of equivalence classes $(A, \lambda, \iota, \eta, \Fil)$
where $(A, \lambda, \iota, \eta) \in \Sh_K(S)$ and $\Fil$ is an (increasing) refinement of
the Hodge filtration given by
\[
  \Fil \, : \, 0 = \mathcal{F}_0 \subseteq \mathcal{F}_1 \subseteq \mathcal{F}_2 \subseteq \mathcal{F}_3 = \Hdr(A)_{\tau_0},
\]
where $\mathcal{F}_r = \omega_{A, \tau_0}$.
Note that $\Fl_K$ is of dimension $3$, regardless of the signature $(r, s)$.
It is instead the dimension of the fibers of $\pi$ that varies with the signature.
\subsubsection{Automorphic line bundles}
The map $\Fl_K \rightarrow \left[* / B\right]$ induces a functor
\[
  \Rep(B) \longrightarrow VB(\Fl_K).
\]
Given a character $\lambda \in X^*(T)$, by first twisting by $w_c$, and then inflating
it to a representation of $B$, we obtain a map
\[
  \mathcal{L} : X^*(T) \xlongrightarrow{w_c} X^*(T) \longrightarrow \Rep(B) \longrightarrow VB(\Fl_K).
\]
The image of this functor takes values in line bundles, since the corresponding representation
of $B$ is a character. The reason for introducing the twist by $w_c$ above is to have
\[
  \pi_* \mathcal{L}(\lambda) = \omega(\lambda),
\]
for dominant $\lambda \in X^*(T)$.
\subsection{The Hodge and conjugate filtrations}
As a preliminary for introducing the \textit{coarse stratification} on our flag Shimura
varieties, let us discuss the Hodge and conjugate filtration of the universal
abelian variety $\mathcal{A}$ in characteristic $p$.
Throughout this subsection we will be working over $\Spec \Fpt$.
\subsubsection{Orthogonal complements}
Given a subbundle $\mathcal{E} \subseteq \Hdr(\mathcal{A})_{\tau_i}$ over $\Fl_K$, let $\mathcal{E}^{\bot} \subseteq \Hdr(\mathcal{A})_{\tau_{i + 1}}$ denote the kernel of the following composition
\[
  \Hdr(\mathcal{A})_{\tau_{i + 1}} \xlongrightarrow{\sim} \Hdr(\mathcal{A})_{\tau_{i + 1}}^{**} \twoheadrightarrow \lambda(\mathcal{E})^* \otimes \delta.
\]
Here we insert a twist by $\delta$ (which is trivial as a line bundle), so that
the orthogonal complement is compatible with the corresponding operation on $B$-subrepresentations of $W_{\tau_0}$ under the morphism
$\Fl_K \rightarrow \left[ * / B \right]$.
Let us start with some basic linear algebra.
\begin{lemma} \label{lem:lin-alg-lemma}
  Given a morphism $\Phi : \mathcal{E} \rightarrow \mathcal{F}$ of vector bundles on some scheme $S$,
  and a subbundle $\mathcal{G} \subseteq \mathcal{F}$. Then we have that
  \[
    \Phi^*(\mathcal{G}^*) = (\Phi^{-1}(\mathcal{G}) / \ker \Phi)^*.
  \]
\end{lemma}
\begin{proof}
  By considering the restriction $\Phi|_{\Phi^{-1}(\mathcal{G})}$, this reduces to
  $\mathcal{G} = \mathcal{F}$. The lemma now follows from the fact that $\Img \Phi^* = (\Img
  \Phi)^*$.
\end{proof}
\begin{corollary} \label{cor:SES-orth-compl}
  For any subbundle $\mathcal{E} \subseteq \Hdr(\mathcal{A})_{\tau_i}$, we can characterize
  $\mathcal{E}^\bot$ as the following kernel
  \[
    0 \longrightarrow \mathcal{E}^\bot \longrightarrow \Hdr(\mathcal{A})_{\tau_{i + 1}} \longrightarrow \lambda^{-1}(\mathcal{E}^*) \longrightarrow 0.
  \]
\end{corollary}
\begin{proof}
  Let us implicitly identify $\Hdr(\mathcal{A})_{\tau_{i + 1}} \cong \Hdr(\mathcal{A})_{\tau_{i + 1}}^{**}$ in this proof.
  The fact that
  $- \lambda = \lambda^*$ and Lemma \ref{lem:lin-alg-lemma} gives us the following chain
  of equalities, from which the corollary follows
  \[
    \lambda(\mathcal{E})^* = - \lambda(\mathcal{E})^* = \lambda^*(\mathcal{E})^* = \lambda^{-1}(\mathcal{E}^*).
  \]
\end{proof}
\begin{corollary}
  For any subbundle $\mathcal{E} \subseteq \Hdr(\mathcal{A})_{\tau_i}$, we have  an isomorphism $\mathcal{E}^{\bot, \bot} \cong \mathcal{E}$ of subbundles in $\Hdr(\mathcal{A})_{\tau_i}$.
\end{corollary}
\begin{proof}
  The above corollary implies that
  $\lambda(\mathcal{E}^{\bot})^* \cong \Hdr(\mathcal{A})_{\tau_i} / \mathcal{E}$, from
  which the desired isomorphism easily follows by unfolding the definition of $\mathcal{E}^{\bot, \bot}$.
\end{proof}
Taking the orthogonal complement of $\Fil$, gives us the following filtration
\[
  \Fil^\bot \, : \, 0 = \mathcal{F}_3^\bot \subseteq \mathcal{F}_2^\bot \subseteq \mathcal{F}_1^\bot \subseteq
  \Hdr(\mathcal{A})_{\tau_1}.
\]
We have the following comparison isomorphisms between the graded pieces of the two filtrations.
\begin{lemma} \label{lem:compl-compare}
  For each $i = 0, 1, 2$, there is an isomorphism of line bundles
  \[
    \mathcal{F}_i^\bot / \mathcal{F}_{i + 1}^\bot \cong (\mathcal{F}_{i + 1}
  / \mathcal{F}_i)^* \otimes \delta.
  \]
\end{lemma}
\begin{proof}
  By untangling the definition, we obtain the following chain of isomorphisms
  \begin{align*}
    \mathcal{F}_i^\bot / \mathcal{F}_{i + 1}^\bot & \cong 
    \ker ( \lambda(\mathcal{F}_{i + 1})^* \otimes \delta \twoheadrightarrow \lambda(\mathcal{F}_{i})^* \otimes \delta ) \\
                                                  & \cong (\lambda(\mathcal{F}_{i + 1}) / \lambda(\mathcal{F}_{i}))^* \otimes \delta \\
                                      & \cong (\mathcal{F}_{i + 1}/\mathcal{F}_{i})^* \otimes \delta.
  \end{align*}
\end{proof}
Together, $\Fil$ and $\Fil^\bot$ assemble to an extension of the Hodge filtration of $\Hdr(A)$.
\subsubsection{The conjugate filtration} \label{sect:conj-filtr}
On $\Sh_K$, we have an \textit{increasing} conjugate filtration given by
\[
  0 \subseteq (\ker V)_{\tau_0} \subseteq \Hdr(\mathcal{A})_{\tau_0}.
\]
Note that it is common to view the conjugate filtration as decreasing, however since we are taking $\tau_0$-isotypic components it is more convenient to view it as an increasing filtration (it makes the comparison in Lemmas \ref{lem:conj-filtr} and \ref{lem:conj-compl} as clean as possible).
There are comparison isomorphisms between the graded pieces of the conjugate filtration,
and the graded pieces of the $(p)$-twisted $\tau_1$-isotypic component of the Hodge filtration.
These are given by
\begin{equation} \label{eq:conj-comparison}
  V : \Hdr(\mathcal{A})_{\tau_0} / (\ker V)_{\tau_0} \xrightarrow{\sim} \omega_{\tau_1}^{(p)} \text{ and } F : (\Hdr(A)_{\tau_1} / \omega_{\tau_1})^{(p)} \xrightarrow{\sim} (\ker V)_{\tau_0}.
\end{equation}
In particular, notice that $(\ker V)_{\tau_0}$ is locally free of rank $r$.
These two filtrations, along with the isomorphisms between the graded pieces, give
$\Hdr(\mathcal{A})$ the structure of a $G$-zip (see §\ref{sect:gzip}).
We now want to extend the conjugate filtration to a full filtration on the flag Shimura variety.
For this purpose, we define
\begin{equation} \label{eq:G-defn}
  \mathcal{G}^i :=
    \begin{cases}
      F(\mathcal{F}_i^{\bot, (p)}) & \text{for } 0 \le i \le r, \\
    V^{-1}(\mathcal{F}_{i}^{\bot, (p)}) & \text{for } r \le i \le r + s.
    \end{cases}
\end{equation}
Note that this is well-defined in the overlapping boundary case, and that it induces \textit{decreasing} filtrations on $(\ker V)_{\tau_0}$ and $\Hdr(\mathcal{A})_{\tau_0} / (\ker V)_{\tau_0}$.
Again, the ordering is chosen to make Lemmas \ref{lem:conj-filtr} and \ref{lem:conj-compl} as clean as possible.
We denote this filtration by $\Filc := \mathcal{G}^\bullet$.

The following lemma describes the graded pieces of $\Filc$ in terms of the graded pieces of $\Fil$.
\begin{lemma} \label{lem:conj-filtr}
  We have isomorphisms
  \[
    \mathcal{G}^i /  \mathcal{G}^{i + 1} \cong (\mathcal{F}_{i + 1} / \mathcal{F}_{i})^{\otimes - p} \otimes \delta^{\otimes p}
  \]
  for all $0 \le i \le r + s$.
\end{lemma}
\begin{proof}
  Since $(\ker F)_{\tau_1} = (\Img V)_{\tau_1} = \omega_{\tau_1}^{(p)}$, this is
  clear from the definition of $\mathcal{G}^\bullet$ and Lemma \ref{lem:compl-compare}.
\end{proof}
We can again take orthogonal complements to obtain the filtration $\Filc^\bot$, which
together with $\Filc$ assembles to a refinement of the conjugate filtration on $\Hdr(A)$.

The following lemma shows that $\Filc^\bot$ can also be obtained analogously to
\eqref{eq:G-defn}.
\begin{lemma} \label{lem:conj-compl}
  For any subbundle $\mathcal{E} \subseteq \Hdr(A)_{\tau_i}$, we have that
  \[
    F(\mathcal{E}^{\bot, (p)}) = V^{-1}(\mathcal{E}^{(p)})^\bot.
  \]
  In particular, we get that
  \[
    \mathcal{G}^{i, \bot} = 
    \begin{cases}
      V^{-1}(\mathcal{F}_{i}^{(p)}) & \text{for } 0 \le i \le r, \\
      F(\mathcal{F}_{i}^{(p)}) & \text{for } r \le i \le r + s.
    \end{cases}
  \]
\end{lemma}
\begin{proof}
  First, Corollary \ref{cor:SES-orth-compl} implies that
  \begin{equation} \label{eq:top-row}
    F(\mathcal{E}^{\bot, (p)}) = \ker\left(F(\Hdr(A)_{\tau_{i + 1}}^{(p)}) \longrightarrow
  F(\lambda^{(p), -1}(\mathcal{E}^{*, (p)}))\right).
  \end{equation}
  Recall that $\lambda \circ F = V^* \circ \lambda^{(p)}$. Hence
  $F(\lambda^{(p), -1}(\mathcal{E}^{*, (p)})) = \lambda^{-1}(V^*(\mathcal{E}^{*, (p)}))$.
  Now Lemma \ref{lem:lin-alg-lemma} implies that $V^*(\mathcal{E}^{*, (p)}) = (V^{-1}(\mathcal{E}^{(p)}) / (\ker V)_{\tau_{i + 1}})^*$.
  Thus we can consider the following commutative diagram, with all rows and columns known to be exact, except the middle row.
  The exactness of the top row follows from \eqref{eq:top-row}, and the exactness of the middle column follows from $(\ker V)_{\tau_i}^\bot = (\ker V)_{\tau_{i + 1}}$.
\[\begin{tikzcd}
	& 0 & 0 & 0 & \\
	0 & {F(\mathcal{E}^{\bot, (p)})} & {(\ker V)_{\tau_i}} & 
  {\lambda^{-1}(V^*(\mathcal{E}^{*, (p)}))} & 0 \\
	0 & {F(\mathcal{E}^{\bot, (p)})} & {\Hdr(A)_{\tau_i}} & {\lambda^{-1}((V^{-1}(\mathcal{E}^{(p)})^*)} & 0 \\
	0 & 0 & {\lambda^{-1}((\ker V)_{\tau_{i + 1}}^*)} & {\lambda^{-1}((\ker V)_{\tau_{i + 1}}^*)} & 0 \\
	&& 0 & 0
	\arrow[from=1-2, to=2-2]
	\arrow[from=1-3, to=2-3]
	\arrow[from=1-4, to=2-4]
	\arrow[from=2-1, to=2-2]
	\arrow[from=2-2, to=2-3]
	\arrow[from=2-2, to=3-2]
	\arrow[from=2-3, to=2-4]
	\arrow[from=2-3, to=3-3]
	\arrow[from=2-4, to=2-5]
	\arrow[from=2-4, to=3-4]
	\arrow[from=3-1, to=3-2]
	\arrow[from=3-2, to=3-3]
	\arrow[from=3-2, to=4-2]
	\arrow[from=3-3, to=3-4]
	\arrow[from=3-3, to=4-3]
	\arrow[from=3-4, to=3-5]
	\arrow[from=3-4, to=4-4]
	\arrow[from=4-1, to=4-2]
	\arrow[from=4-2, to=4-3]
	\arrow[from=4-3, to=4-4]
	\arrow[from=4-3, to=5-3]
	\arrow[from=4-4, to=4-5]
	\arrow[from=4-4, to=5-4]
\end{tikzcd}\]
The middle row is exact by the $9$-lemma, from which $F(\mathcal{E}^{\bot, (p)}) = V^{-1}(\mathcal{E}^{(p)})^\bot$ follows.
\end{proof}

\subsection{The coarse stratification}
In characteristic $p$, flag Shimura varieties admit a well-behaved stratification known as the
\textit{coarse stratification} \cite{flag-strata}.
It generalizes the classical $a$-number stratification of a Shimura variety.

In this subsection, let us give a brief indication of the general theory of this stratification, and
its relation to the stack of $G$-zipflags.
Let us remark that once we specialize to the type $(2, 1)$ and $(3, 0)$-case,
we will be able to give concrete descriptions of these strata (and the partial Hasse invariants), independent of the
theory of $G$-zips.
Hence, a deep understanding of stacks of $G$-zips will not be necessary for reading this paper.

Throughout this subsection we continue to work over $\Spec \Fpt$.
\subsubsection{The stack of $G$-zips} \label{sect:gzip}
Consider the algebraic group
\[
  E := \left\{ (x, y) \in P \times P^{-, (p)} \mid \varphi(\overline{x}) = \overline{y} \right\}.
\]
Here  $\overline{x}$ denotes the projection of $x$ to the Levi $M \subseteq P$ (resp. to $M^{(p)} \subseteq P^{-, (p)}$), and $\varphi : M \to M^{(p)}$ denotes the Frobenius homomorphism.
The group $E$ acts on $G$ via $(x, y) \cdot g := xgy^{-1}$. From this, we can define the stack of $G$-zips as
\[
  \Gzip^\mu := [E \backslash G].
\]
This is a smooth algebraic stack of dimension $0$, and
it parametrizes objects known as $G$-zips \cite[Theorem 1.5]{Gzip}.
Informally, a $G$-zip should be thought of as a
$G$-torsor $V$, equipped with two filtrations $\mathcal{F}_\bullet$ and $ \mathcal{G}^\bullet$ (whose types are determined by $P$ resp. $P^-$), as
well as a Frobenius semi-linear isomorphism between their associated graded objects.
We have already seen an example of such an object: $\Hdr(\mathcal{A})$ equipped
with the Hodge and the conjugate filtration, as well as the isomorphisms from \eqref{eq:conj-comparison}.
This is the only example we will care about.
It follows from \cite{zhang} that this induces a
smooth morphism $\zeta : \ShKtwo \rightarrow \Gzip^\mu$, through which
$\ShKtwo \rightarrow \left[ * / P \right]$ factors.
In particular, we can obtain the automorphic vector bundles via pullback from $\Gzip^\mu$.

We will also need an analogue of $\Gzip^\mu$ for the flag Shimura varieties.
Given a parabolic $Q \subseteq P$, we can define $\Gzf^{\mu, Q}$ via pullback
\[\begin{tikzcd}
	{\Gzf^{\mu, Q}} & {\left[* / Q\right]} \\
	{\Gzip^\mu} & {\left[* / P\right].}
	\arrow[from=1-1, to=1-2]
	\arrow[from=1-1, to=2-1]
	\arrow[from=1-2, to=2-2]
	\arrow[from=2-1, to=2-2]
\end{tikzcd}\]
This stack parametrizes so-called $G$-zipflags, which by definition are $G$-zips
equipped with a refinement of the Hodge filtration.
As usual, we will omit $Q$ from the notation when $Q = B$.
As an immediate consequence of the definition,
we also get the following pullback square
\[\begin{tikzcd}
  {\Fl_{K, Q}} & {\Gzf^{\mu, Q}} \\
	{\Sh_K} & {\Gzip^\mu.}
	\arrow[from=1-1, to=1-2]
	\arrow["\pi", from=1-1, to=2-1]
	\arrow[from=1-2, to=2-2]
	\arrow[from=2-1, to=2-2]
\end{tikzcd}\]
It follows as in §\ref{sect:conj-filtr}, that
given a refinement of the Hodge filtration,
we can define a suitable refinement of the conjugate filtration.
Given a $G$-zipflag, we will abusively let $\mathcal{F}_\bullet$ (or $\Fil$) denote
its refined Hodge filtration (i.e. the one of type $Q$), and we let $\mathcal{G}^\bullet$ (or $\Filc$)
denote the induced refinement of the conjugate filtration.
\subsubsection{The coarse stratification}
The coarse stratification on $\Gzf^{\mu, Q}$ (and hence on $\Fl_{K, Q}$),
is given by measuring the relative position between $\mathcal{F}_\bullet$ and $\mathcal{G}^\bullet$.
Let us explain this only in the case $Q = B$ for simplicity.
The stratification can be obtained geometrically from the following morphism of stacks
\[
  \Gzf^\mu \longrightarrow [B \backslash G / B],
\]
given by sending a $G$-zipflag to the relative position of $\mathcal{F}_\bullet$ and $\mathcal{G}^\bullet \cdot z$, where $z := w_0 w_c$ (viewed as an element in $G$).
The point of $z$ is to conjugate $\mathcal{G}^\bullet$ to have the same stabilizer as
$\mathcal{F}_\bullet$.
The strata on $\Gzf^\mu$ are then given as the preimages of points in $[B \backslash G / B]$, which
by the Bruhat decomposition correspond (uniquely) to double cosets $BwB$ for $w \in W(T)$.
It is well-known that $\overline{BwB} = \cup_{w' \le w} Bw'B$, where the order is given by the
Bruhat order on $W(T)$. Thus, the closure relations between the strata are exactly determined by the
Bruhat order on $W(T)$.

Given $w \in W(T)$ we will let $\mathcal{C}_w$ denote the corresponding stratum in $\Gzf^\mu$ (i.e. the preimage of $BwB$), and we will let $\overline{\mathcal{C}}_w$ denote its closure.
Let $\mathcal{Z}_w \subseteq \Fl_K$ denote the preimage of $\mathcal{C}_w$, and we let $\Zcb_w$ denote its closure.
From \cite[Proposition 6.3.1]{flag-strata} we have that $\mathcal{Z}_w$ has dimension $\ell(w)$.
It follows from the explicit description of $W(T) \cong S_3$, together with its Bruhat order, that our flag Shimura variety $\FlKtwo$ will have $6$ strata: one $3$-dimensional stratum, two $2$-dimensional, two $1$-dimensional, and one $0$-dimensional stratum. Each closed stratum of a given dimension
contains all strata of strictly smaller dimension.

\subsubsection{Partial Hasse invariants}
The coarse stratification satisfies the property that for each $w \in W$, each codimension $1$ stratum $\overline{\mathcal{Z}}_{w'} \subseteq \overline{\mathcal{Z}}_w$ is given by the vanishing of a section of some line bundle over $\mathcal{\overline{Z}}_w$.
Such sections are known as \textit{partial Hasse invariants}.
We will be able to explicitly write
these down whenever we need to use them, so we will not go much into the general theory.

\subsection{The Shimura variety of type $(2, 1)$} \label{sect:2-1}
Let us now specialize the above theory to the case $(r, s) = (2, 1)$.
From now on, we will reuse the notation from the general type $(r, s)$ setting
in the type $(2, 1)$-case, and we will add a prime $'$ to denote corresponding object in the type $(3, 0)$-case (which will be discussed in §\ref{sect:3-0}).
In particular, we will let
$(E, V, \langle \cdot, \cdot \rangle, h, \mathcal{O}_{E_p}, \Lambda)$ denote
the integral Shimura datum in the type $(2, 1)$-case, and $G$ will denote the associated group.
Fix a small enough level structure $K = K^pK_p \subseteq G(\adf)$
that is hyperspecial at $p$.

For the remainder of this paper,
we will be working in characteristic $p$.
Therefore, let us use $\ShKtwo$ to denote the special fiber of the integral model
of level $K$. We let $\mathcal{A}$ denote the associated universal abelian variety,
and we write $\FlKtwo$ for the corresponding flag Shimura variety (over $\Spec \Fpt$).
\subsubsection{Some group theory}
In the type $(2, 1)$-case our Hodge cocharacter $\mu$ is given by
$z \mapsto \operatorname{diag}(z, z, 1) \times z$.
Thus, the Levi $M$ is given by the subgroup of matrices of the form
\[
  \begin{pmatrix}
  * & * & 0\\
  * & * & 0\\
  0 & 0 & *\\
  \end{pmatrix} \times (*),
\]
and $P$ will consist of matrices of the form
\[
  \begin{pmatrix}
  * & * & *\\
  * & * & *\\
  0 & 0 & *\\
  \end{pmatrix} \times (*).
\]
\subsubsection{The filtrations}
In the type $(2, 1)$-case, we will write the Hodge filtration $\Fil$ on $\FlKtwo$ as
\[
  0 \subseteq \mathcal{F} \subseteq \omega_{\mathcal{A}, \tau_0} \subseteq \Hdr(\mathcal{A})_{\tau_0}.
\]
We will reserve the notation $\mathcal{F}_i$ for the type $(3, 0)$-case (the above filtration
will not be used after this subsection, so we hope that no confusion will occur).
The conjugate filtration $\Filc$ is then given by
\begin{equation} \label{eq:conj-filtr-21-def}
  0 \subseteq \mathcal{G} := F(\mathcal{F}^{\bot, (p)}) \subseteq (\ker V)_{\tau_0} \subseteq
    \Hdr(\mathcal{A})_{\tau_0}.
\end{equation}
\subsubsection{Automorphic vector bundles}
The parabolic $P$ induces a $2$-dimensional $\Fpt$-subspace $F \subseteq W_{\tau_0}$.
This allows us to identify $M \cong \GL(F) \times \GL(W_{\tau_0} / F) \times
\mathbb{G}_{m, \Fpt}$. With this identification, we get the following description for
dominant $\lambda = (a, b, c, d)$ (i.e. $a \ge b$)
\[
  \Ind_{M \cap B}^{M} w_c\lambda \cong 
  \Sym^{a - b} F \otimes (\det F)^{\otimes b} \otimes (W_{\tau_0} / F)^{\otimes c} \otimes \delta^{\otimes d}.
\]
It then follows that (for $\lambda$ dominant)
\[
  \omega(\lambda) \cong
  \Sym^{a - b} \omega_{\mathcal{A}, \tau_0} \otimes (\det \omega_{\mathcal{A}, \tau_0})^{\otimes b} \otimes (\Lie \mathcal{A}^\vee)_{\tau_0}^{\otimes c} \otimes \delta^{\otimes d}.
\]
\subsubsection{Automorphic line bundles}
Tracing through the construction of the automorphic line bundles given above (and being
mindful of the twist by $w_c$), we get the following description in terms
of the graded pieces of $\Fil$
\begin{equation} \label{eq:lb-defn}
  \mathcal{L}(a, b, c, d) \cong (\omega_{\mathcal{A}, \tau_0} / \mathcal{F})^{\otimes a} \otimes \mathcal{F}^{\otimes b} \otimes (\Lie \mathcal{A}^\vee)_{\tau_0}^{\otimes c} \otimes \delta^{\otimes d},
\end{equation}
for any dominant $(a, b, c, d) \in X^*(T)$.

\subsubsection{Closed strata and partial Hasse invariants} \label{sect:2-1-hasses}
We have two closed strata of codimension $1$, given by the conditions $\mathcal{F} \subseteq (\ker
V)_{\tau_0}$ resp. $\mathcal{G} \subseteq \omega_{\mathcal{A}, \tau_0}$.
One can check that these correspond to the elements $s_2s_1 = (132)$ resp. $s_1s_2 = (123)$ of $W(T) = S_3$, so we will
label these by $\Zcb_{21}$
resp. $\Zcb_{12}$.
Notice that these are given by the vanishing loci of the projection morphisms
\[
  \ha_{21} : \mathcal{F} \longrightarrow \Hdr(\mathcal{A})_{\tau_0} / (\ker V)_{\tau_0},
\]
resp.
\[
  \ha_{12} : \mathcal{G} \longrightarrow \Hdr(\mathcal{A})_{\tau_0} / \omega_{\mathcal{A}, \tau_0}. 
\]
Note that after re-arranging (and applying Lemma \ref{lem:conj-filtr}),
we can consider these as sections of automorphic line bundles: $\ha_{21} \in H^0(\Fl_K, \mathcal{L}(0,-1, -p, p))$ and $\ha_{12} \in H^0(\Fl_K, \mathcal{L}(p, 0, 1, - p))$.

There are two closed strata in $\FlKtwo$ of codimension $2$. They are defined by the closed conditions 
$\omega_{\mathcal{A}, \tau_0} \subseteq (\ker V)_{\tau_0}$ resp. $\mathcal{F} \subseteq \mathcal{G}$.
Again, one may check that these correspond to the elements $s_1$ resp. $s_2$ in $W(T)$,
so we will denote them by $\Zcb_{1}$ resp. $\Zcb_{2}$.
Notice that, since these are containments of subbundles of the same rank, they are actually equalities.
Considering both of these strata as subschemes of $\Zcb_{12}$, we get
corresponding partial Hasse invariants
\[
  \ha_{1} : (\ker V)_{\tau_0} / \mathcal{G} \longrightarrow \Hdr(\mathcal{A})_{\tau_0} / \omega_{\mathcal{A}, \tau_0}
\]
and
\[
  \ha_{2} : \mathcal{G} \longrightarrow \omega_{\mathcal{A}, \tau_0} / \mathcal{F}.
\]
We can also consider these as sections of 
$\mathcal{L}(0,p,1,-p)|_{\Zcb_{12}}$ resp. $\mathcal{L}(p + 1, 0, 0, -p)|_{\Zcb_{12}}$.
\begin{lemma} \label{lem:vanish-simply}
  We have $\operatorname{div}(\ha_i) = - \Zcb_i$ for $i \in \{21, 12, 1, 2\}$,
  where we view $\ha_i$ as a section over $\Zcb_{12}$ when $i \in \{ 1, 2 \}$.
\end{lemma}
\begin{proof}
It is clear that these sections can be obtained via pullback from
$[ B \backslash G / B ]$
  (their definition only depends on the two filtrations and their relative positions).
  Next, since the morphism $\FlKtwo \rightarrow [ B \backslash G / B ]$ is smooth,
  it follows that it suffices to compute the corresponding divisors on $[ B \backslash G / B ]$.
  This can be done by combining Theorem 2.2.1c) and Lemma 3.1.1b) in \cite{strata-hasses}.
\end{proof}
We will say that an abelian variety is superspecial if it is isomorphic to a
product of supersingular elliptic curves.
\begin{lemma} \label{lem:P1-fibers}
  The stratum $\Zcb_{1}$ is a $\mathbb{P}^1$-bundle over the $0$-dimensional superspecial locus of $\ShKp$.
\end{lemma}
\begin{proof}
  Notice that the condition $\omega_{\mathcal{A}, \tau_0} \subseteq (\ker V)_{\tau_0}$
is independent of the subbundle $\mathcal{F}$. Hence, we already get that
  $\Zcb_{1}$ must be a $\mathbb{P}^1$-bundle over some locus in the base $\ShKp$.
  Next, it follows by taking orthogonal complements that the condition defining $\Zcb_{1}$ is equivalent
  to demanding $\omega_{A} = \ker V$.
  It follows from \cite[Theorem 2]{Oort1975} that this condition is equivalent to the abelian variety being
  superspecial.
\end{proof}

\subsection{The Shimura set of type $(3, 0)$} \label{sect:3-0}
Let $(E, V, \langle \cdot, \cdot \rangle', h', \mathcal{O}_{E_p}, \Lambda)$ denote
the integral PEL datum in the type $(3, 0)$ case, and let $G'$ denote the associated group.
Since we are working in odd dimension, there is only one unitary group over a local non-archimedean field.
In particular, we have that $G_{\Q_\ell} \cong G'_{\Q_\ell}$ for all primes $\ell$.
Thus, we can pick a level $K' \subseteq G'(\adf) \cong G(\adf)$ corresponding to $K$.
Let $\ShKtr$ denote the special fiber of the corresponding Shimura set of level $K'$.
Since $G'_{\Qp} \cong G_{\Qp}$, we will generally replace $G'_{\Qp}$ (and its special fiber) with
$G_{\Qp}$ (resp. $G_{\Fpt})$) in what follows. In particular, we can use $X^*(T)$ as the weights
of automorphic vector bundles in both settings.
Moreover, in the type $(3, 0)$-case, note that the stabilizer of the Hodge cocharacter (over $\Spec \Fpt$) is given by $P' = G_{\Fpt}$, and the
corresponding Levi is then also $G_{\Fpt}$.
We let $\FlKtr$ denote the special fiber of the corresponding flag Shimura
variety, and $\mathcal{A}'$ denote the universal abelian variety.

\subsubsection{Abelian varieties of type $(3, 0)$}
Fix an abelian variety $A$ of type $(3, 0)$ over a perfect field $k$ in characteristic $p$.
The signature condition implies in particular that $\omega_A = \omega_{A, \tau_0} = (\ker V)_{\tau_0} = \ker V$.
Thus, $A$ must be superspecial.

It follows from \cite[Proposition 2.1.2.]{gu2n2} that there exists some supersingular
elliptic curve $\mathcal{E}$, and some self-dual $\mathcal{O}_{E_p}$-Hermitian lattice $\Lambda_0$ such that
$\mathcal{D}(A) \cong \mathcal{D}(\mathcal{E}) \otimes_{\mathcal{O}_{E_p}} \Lambda_0$.
In particular, we deduce from the case of supersingular elliptic curves that
$F^{-1, (p)} \circ V = \sigma^2 : \mathcal{D}(A) \rightarrow \mathcal{D}(A)^{(p^2)}$
and $V^{-1, (p)} \circ F = \sigma^{-2} : \mathcal{D}(A)^{(p^2)} \rightarrow \mathcal{D}(A)$.

\subsubsection{Automorphic vector bundles} \label{sect:aut-vb-threezero}
In the type $(3, 0)$-setting, the automorphic vector bundles are obtained via pullback from $\ShKtr \rightarrow [* / G_{\Fpt} ]$.
For a representation $\sigma \in \Rep(G_{\Fpt})$, we let $\omega'(\sigma)$ denote
the associated automorphic vector bundle.
For a character $\lambda \in X^*(T)$ we let $\omega'(\lambda)$ denote the automorphic vector bundle of
coherent weight $\lambda \in X^*(T)$ as before.
Notice that $\omega'(\lambda) = \omega'(W(\lambda))$, where $W(\lambda)$ denotes the dual Weyl
module associated to $\lambda$.
By using the description of the Shimura set $\ShKtr$ as an adelic double quotient, one may show that for any $\sigma \in \Rep(G_{\Fpt})$, we have
\[
  H^0(\omega'(\sigma), \ShKtr) = \{ f : G'(\Q) \backslash G'(\adf) / K'^p \rightarrow \sigma : f(gk) = k^{-1}f(g) \quad \forall g \in G'(\adf), k \in K'_p \}.
\]
Thus, global sections of automorphic vector bundles for Shimura sets are just spaces of algebraic automorphic forms in the sense of \cite{gross}.

\subsubsection{Automorphic line bundles} We get the following description of the automorphic line bundles
\begin{equation}
  \mathcal{L}'(a, b, c, d) \cong (\mathcal{F}_3 / \mathcal{F}_2)^{\otimes a} \otimes (\mathcal{F}_2 / \mathcal{F}_1)^{\otimes b}
  \otimes \mathcal{F}_1^{\otimes c} \otimes {\delta}^{\otimes d},
\end{equation}
for any $(a, b, c, d) \in X^*(T)$.

\subsubsection{The filtrations}
In the type $(3, 0)$-case, we have the following Hodge and conjugate filtrations.
\[
  \Fil \, : \, 0 \subseteq \mathcal{F}_1 \subseteq \mathcal{F}_2 \subseteq \mathcal{F}_3 = \omega_{\mathcal{A}, \tau_0} = \Hdr(\mathcal{A}')_{\tau_0}
\]
and
\[
  \Filc \, : \, 0 \subseteq F(\mathcal{F}_2^\bot) \subseteq F(\mathcal{F}_1^\bot) \subseteq (\ker V)_{\tau_0} = \Hdr(\mathcal{A}')_{\tau_0}.
\]
We will use the notation $\mathcal{G}^i = F(\mathcal{F}_i^{\bot, (p)})$ as before.
The following lemma shows that,
in the type $(3, 0)$-case, we can completely describe the action of $F$ and $V$ on
the two filtrations
\begin{lemma} \label{lem:frob-ver-reln}
  We have that
  \begin{enumerate}[(i)]
    \item $V(\mathcal{G}^{i, \bot}) = \mathcal{F}_i^{(p)}$
    \item $F(\mathcal{G}^{i, \bot, (p)}) = \mathcal{F}_i$
    \item $F(\mathcal{F}_i^{\bot, (p)}) = \mathcal{G}^{i}$
    \item $V(\mathcal{F}_i^{\bot}) = \mathcal{G}^{i, (p)}$
  \end{enumerate}
\end{lemma}
\begin{proof}
  (i) follows from Lemma \ref{lem:conj-compl}, and (iii) is true by definition.
  Let us now prove (ii). Recall that $V^{-1} \circ F^{(p)} = \sigma^{-2} : \Hdr(\mathcal{A}')_{\tau_1}^{(p^2)} \rightarrow \Hdr(\mathcal{A}')_{\tau_1}$.
  Thus, $V^{-1} \circ F^{(p)}(\mathcal{G}^{i, \bot, (p^2)}) = \mathcal{G}^{i, \bot}$.
  Applying $V$ to both sides, combined with (i), gives
  \[
    (F(\mathcal{G}^{i, \bot, (p)}))^{(p)} = F^{(p)}(\mathcal{G}^{i, \bot, (p^2)}) = V(\mathcal{G}^{i, \bot}) = \mathcal{F}_i^{(p)}.
  \]
  Thus, by untwisting, we obtain (ii). Finally, (iv) follows from the same argument using (iii) instead of (i) (but this time we do not
  even need to untwist).
\end{proof}
\subsubsection{The stratification}
We will mainly be interested in two of the closed strata in $\FlKtr$.
The first stratum, which we will denote by $\Zcb'_{21} \subset \FlKtr$,
is defined by the codimension $1$ condition
\begin{equation} \label{eq:z21'-cond}
  \mathcal{F}_1 \subseteq \mathcal{G}^1.
\end{equation}
It corresponds to the element $s_2s_1 = (132)$ in $W(T)$.
Next, we consider the codimension $2$ stratum given by
\[
  \mathcal{F}_2 = \mathcal{G}^1.
\]
This corresponds to the element $s_1 \in W(T)$, so we denote it by $\Zcb'_1$.
Note that it is contained in $\Zcb'_{21}$ (of codimension $1$).
Notice that the equality $\mathcal{F}_2 = \mathcal{G}^1$ over $\Zcb'_1$ gives us an isomorphism $\mathcal{G}^0 / \mathcal{G}^1 \xrightarrow{\sim} \mathcal{F}_3 / \mathcal{F}_2$.
After tensoring by $(\mathcal{G}^0 / \mathcal{G}^1)^{\otimes -1}$ on both sides,
we obtain an invertible section 
\begin{equation} \label{eq:invertible-hasse}
  \mathcal{O}_{\Zcb'_1} \xrightarrow{\sim} \mathcal{L}'(1, 0, p, -p)|_{\Zcb'_1}.
\end{equation}
This can be viewed as an ``invertible partial Hasse invariant'', and will allow
us to shift the weight of an automorphic line bundle by multiples of $(1, 0, p, -p)$ without affecting the isomorphism
class over $\Zcb'_1$.

Let us also write $\overline{\mathcal{C}}'_i$ to denote the closed stratum of $\Gzf^{\mu'}$
corresponding to $\Zcb'_i$, and $\mathcal{C}'_i$ for the unique open stratum contained in $\overline{\mathcal{C}}'_i$.

\subsubsection{Relation to Deligne--Lusztig varieties}
There is a close relation between the coarse stratification on $\FlKtr$
and the Deligne--Lusztig stratification on $G / B$, as we now explain.
Recall the $\Fpt$-Hermitian space $W = \Lambda / p \Lambda$.
The above discussion on stratifications then applies equally to the flag variety $G / B = \GU(W)_{\Fpt} / B$, after replacing $\Hdr(\mathcal{A}')$ with $W$, and replacing the Frobenius morphism
$F_{\mathcal{A}}$ acting on the Hodge filtration $\mathcal{F}_\bullet$ by the Frobenius of $\GU(W)_{\Fpt}$
acting on the quotient $\GU(W)_{\Fpt} / B$. In particular, we obtain a morhpism $G / B \rightarrow \Gzf^{\mu'}$. Let $X_w$ denote
the preimage of $\mathcal{C}'_w$. This gives a stratification of $G / B$, and
the strata $X_w$ for $w \in W(T)$ are exactly the Deligne--Lusztig varieties.
We will let $\overline{X}_i$ denote the closed stratum corresponding to $\overline{\mathcal{C}}'_i$.
\begin{remark}
Since $\FlKtr$ is just a disjoint union of flag varieties $G / B$, it
follows that $\Zcb'_i$ is just a disjoint union of closed Deligne--Lusztig varieties isomorphic to
$\overline{X}_i$.
\end{remark}

\subsubsection{The partial flag variety} \label{sect:partial-flag}
An important result in this paper is the computation of global sections of certain
automorphic line bundles on $\Zcb'_1$.
For this computation, it will be more convenient to realize $\Zcb'_1$ as a codimension $1$
stratum inside a partial flag Shimura variety (as opposed to a codimension $2$ stratum in the full flag
Shimura variety). Let $P^* \subseteq G_{\Fpt}$ denote the parabolic consisting of matrices of the form
\[
  \begin{pmatrix}
  * & * & *\\
  0 & * & *\\
  0 & * & *\\
  \end{pmatrix} \times (*).
\]
The Hodge filtration on the corresponding flag Shimura variety $\Fl_{K', P^*}$ is
then given by
\[
  \Fil \, : \, 0 \subseteq \mathcal{F}_1 \subseteq \mathcal{F}_3 = \Hdr(\mathcal{A}')_{\tau_0}.
\]
If we let $\mathbb{P}^{2, *}$ denote the projective space parametrizing $2$-dimensional quotients (we reserve $\mathbb{P}^2$ to denote the space parametrizing $1$-dimensional quotients, as is more standard),
then $\Fl_{K', P^*}$ is just given by a disjoint union of varieties isomorphism to $\mathbb{P}^{2, *}$.

There is a projection map $\pi_{B, P^*} : \FlKtr \rightarrow \Fl_{K', P^*}$ given by sending
$\mathcal{F}_\bullet$ to $\mathcal{F}_1$. By abuse of notation, this restricts
to a closed immersion $\pi_{B, P^*} : \Zcb'_1 \rightarrow \Fl_{K', P^*}$.
Indeed, the map is proper, and
since $\mathcal{F}_2 = \mathcal{G}^1$ on $\Zcb'_1$, we see that $\mathcal{F}_2$ is
uniquely determined from $\mathcal{F}_1$, so the map is also a monomorphism.
Hence, we can (and will) treat $\Zcb'_1$ as a closed codimension $1$ stratum inside $\Fl_{K', P^*}$.
\begin{remark} \label{rmk:DL-var}
  It follows from the same argument that the corresponding closed Deligne--Lusztig
  variety $\overline{X}_1$ can be identified with a closed subvariety of $G / P^*$.
  This is exactly (at least after replacing $P^*$ by $P$) the generalized Deligne--Lusztig
  variety considered in \cite{vollaard}.
\end{remark}
We can also consider automorphic vector bundles on $\Fl_{K', P^*}$. In particular,
we will be interested in the line bundles
\[
  \mathcal{L}'_{P^*}(a, a, b, d) := \det(\mathcal{F}_3 / \mathcal{F}_1)^{\otimes a} \otimes
  \mathcal{F}_1^{\otimes b} \otimes \delta^{\otimes d}.
\]
Notice that $\pi_{B, P^*}^*\mathcal{L}'_{P^*}(a, a, b, d) \cong \mathcal{L}'(a, a, b, d)$.

\section{Geometric Jacquet--Langlands} \label{sect:exotic-hecke-main}
In this section we describe the geometric Jacquet--Langlands phenomenon,
relating the unitary Shimura variety of type $(2, 1)$ to the unitary Shimura set of type $(3, 0)$,
by identifying the two $1$-dimensional closed strata $\Zcb_{2}$ and $\Zcb'_1$.

\subsection{The family of isogenies over $\Zcb'_1$} \label{sect:exotic-hecke}
Recall that $\mathcal{A}'$ denotes the universal abelian variety over $\FlKtr$.
Our goal is to define a family of isogenies $\mathcal{A}'|_{\Zcb'_{21}} \rightarrow \mathcal{A}_3$, such that $\mathcal{A}_3$ obtains a prime-to-$p$ quasi-polarization 
and an $\mathcal{O}_E$-action of signature $(2, 1)$.
After restricting to $\Zcb'_1$, this family of isogenies
will allow us to identify $\Zcb'_1$ with $\Zcb_{2}$.

The following lemma will allow us to construct subgroups from
quotients of the co-Lie algebra of an abelian variety.
\begin{lemma} \label{lem:integrate-subgroup}
  Given a scheme $S$ of characteristic $p$, there exists a contravariant equivalence of
  categories
  \[
    \left\{ 
      \begin{array}{l}
        \text{p-torsion finite flat group} \\
        \text{schemes } G/S \text{ with } F_G = 0
      \end{array}
    \right\}
    \longleftrightarrow
    \left\{
    (\mathcal{M}, V : \mathcal{M} \rightarrow \mathcal{M}^{(p)}) \, | \, \mathcal{\M} \text{ locally free over } S
    \right\}.
  \]
  The functor is given by sending $G \mapsto (\omega_G, V_G : \omega_G \rightarrow \omega_G^{(p)})$.
\end{lemma}
\begin{proof}
  This is \cite[Example 2.3.]{dejong}.
\end{proof}
The following observation will also be useful
\begin{lemma} \label{lem:ht-one}
  Given a scheme $S/\Spec \Fpt$, and a finite flat group scheme $G/S$ of order $p$. Then
  if $G$ has an $\mathcal{O}_{E}$-action it must be an $F, V$-kernel.
\end{lemma}
\begin{proof}
  We can check this pointwise, so it suffices to assume that $S = \Spec k$ for some algebraically closed field $k/\Fpt$. Next, note that
  $\mathcal{D}(G)_{\Spec k}$ is $1$-dimensional over $k$. Hence, it must be isomorphic to an isotypic
  component under the $\mathcal{O}_{E}$-action. Since $F_G$ and $V_G$ exchanges
  isotypic components, they must both be zero.
\end{proof}
\begin{remark}
  This lemma shows in particular that $\ShKtwo$ has empty ordinary locus.
\end{remark}
The construction will now proceed in three steps. First, we will construct a sequence of order $p$ isogenies of abelian varieties over $\FlKtr$
\[
  \mathcal{A}' \longrightarrow \mathcal{A}_1 \longrightarrow \mathcal{A}_2.
\]
After restricting to the stratum $\Zcb'_{21}\subseteq \FlKtr$, we can define the final order $p$ isogeny
\[
  \mathcal{A}_2|_{\Zcb'_{21}} \longrightarrow \mathcal{A}_3.
\]

The reason for breaking up the isogeny $\mathcal{A}' \rightarrow \mathcal{A}_3$ into pieces,
is to ensure that the kernel at each step is an $F, V$-kernel (this will not be true for the composition).
This will allow us to apply the lemmas from §\ref{sect:ht-one}.
Let us remark that the kernel of the isogeny $\mathcal{A}' \rightarrow \mathcal{A}_2$ \textit{will} be an $F,
V$-kernel,
so in principle we only need to factor $\mathcal{A}' \rightarrow \mathcal{A}_3$ into a composition of two isogenies.
However, we still choose to break it up into three pieces, as it illustrates that the isogeny
$\mathcal{A}_1 \rightarrow \mathcal{A}_3$ gives a $\Zcb'_{21}$-valued point of the Iwahori level Shimura
variety of type $(2, 1)$.

An important point of this construction is that it should allow us to identify
$\mathcal{L}(\lambda)|_{\Zcb_{2}}$, for $\lambda \in X^*(T)$, with $\mathcal{L}'(\lambda')|_{\Zcb'_1}$, for some explicit $\lambda' \in X^*(T)$.
Thus, it will be important to at each step relate the isotypic decomposition of
$\omega_{\mathcal{A}_i}$ and $\Lie (\mathcal{A}_i^\vee)$ to automorphic line bundles
restricted from $\FlKtr$.
\subsubsection{The first isogeny} \label{sect:first-quot}
Consider the surjection
\[
  \omega_{\mathcal{A}'} \twoheadrightarrow \mathcal{F}_3 / \mathcal{F}_2.
\]
Since $\omega_{\mathcal{A}'} = \ker V_{\mathcal{A}'}$,
we can apply Lemma \ref{lem:integrate-subgroup} to obtain an $F, V$-kernel
$G_0 \hookrightarrow \mathcal{A}'$. Let $\mathcal{A}_1$ denote $\mathcal{A}' / G_0$.
Note that it immediately inherits an $\mathcal{O}_{E}$-action, since $G_0$ is $\mathcal{O}_{E}$-stable.
Applying Lemma \ref{lem:hodge-seq-final}, we get the following diagram with
exact rows, describing the Hodge filtrations
\begin{equation} \label{eq:hodge-filtr-1}
\begin{tikzcd}
  0 & {(\mathcal{F}_3 / \mathcal{F}_2)^{(p)}} & {\omega_{\mathcal{A}_1}} & {\omega_{\mathcal{A}'}} & {\mathcal{F}_3/\mathcal{F}_2} & 0 \\
  0 & {\mathcal{H}^0(\Delta(G_0)_{\FlKtr})} & {\Hdr(\mathcal{A}_1)} & {\Hdr(\mathcal{A}')} & {\mathcal{D}(G_0)_{\FlKtr}} & 0 \\
  0 & {\mathcal{F}_3 / \mathcal{F}_2} & {\Lie(\mathcal{A}_1^\vee)} & {\Lie(\mathcal{A}'^\vee)} & {(\mathcal{F}_3 / \mathcal{F}_2)^{(p)}} & {0.}
	\arrow[from=1-1, to=1-2]
	\arrow[from=1-2, to=1-3]
	\arrow["\simeq", from=1-2, to=2-2]
	\arrow[from=1-3, to=1-4]
	\arrow[from=1-3, to=2-3]
	\arrow[from=1-4, to=1-5]
	\arrow[from=1-4, to=2-4]
	\arrow[from=1-5, to=1-6]
	\arrow["0", from=1-5, to=2-5]
	\arrow[from=2-1, to=2-2]
	\arrow[from=2-2, to=2-3]
	\arrow["0", from=2-2, to=3-2]
	\arrow[from=2-3, to=2-4]
	\arrow[from=2-3, to=3-3]
	\arrow[from=2-4, to=2-5]
	\arrow[from=2-4, to=3-4]
	\arrow[from=2-5, to=2-6]
	\arrow["\simeq", from=2-5, to=3-5]
	\arrow[from=3-1, to=3-2]
	\arrow[from=3-2, to=3-3]
	\arrow[from=3-3, to=3-4]
	\arrow[from=3-4, to=3-5]
	\arrow[from=3-5, to=3-6]
\end{tikzcd}
\end{equation}
Furthermore, we also get that the bottom right vertical isomorphism is
induced by $V_{\mathcal{A}'}$.
Hence, by Lemma \ref{lem:frob-ver-reln} i), it is natural to identify 
the quotient
$\Hdr(\mathcal{A}') \rightarrow \mathcal{D}(G_0)_{\FlKtr}$
with the quotient $\Hdr(\mathcal{A}') \rightarrow \mathcal{G}^{3, \bot} / \mathcal{G}^{2, \bot}$.
We will also identify $\mathcal{H}^0(\Delta(G_0)_{\FlKtr}) \cong \mathcal{G}^{3, \bot} / \mathcal{G}^{2, \bot}$ (using Lemma \ref{lem:frob-ver-descr}).
Note that all arrows in \eqref{eq:hodge-filtr-1} are $\mathcal{O}_{E}$-equivariant,
hence we can deduce the following decompositions into isotypic components
\begin{equation} \label{eq:decomp-first-vb}
  \omega_{\mathcal{A}_1} \cong \omega_{\mathcal{A}_1, \tau_0} \oplus \omega_{\mathcal{A}_1, \tau_1} \cong \mathcal{F}_2 \oplus \mathcal{G}^{3, \bot} / \mathcal{G}^{2, \bot}
\end{equation}
and
\[
  \Lie(\mathcal{A}_1^\vee) \cong \Lie(\mathcal{A}_1^\vee)_{\tau_0} \oplus \Lie(\mathcal{A}_1^\vee)_{\tau_1} \cong \mathcal{F}_3 / \mathcal{F}_2 \oplus \mathcal{G}^{2, \bot}.
\]
In particular, $\mathcal{A}_1$ is of type $(2, 1)$.
\begin{remark} \label{rmk:special}
  $\mathcal{A}_1$ can be equipped with a quasi-polarization making it of
  \textit{special level} at $p$ (i.e. its level at $p$ corresponds to the other, non-hyperspecial, maximal parahoric). Since the construction only depends on 
  $\mathcal{F}_2$ and not all of $\mathcal{F}_\bullet$, this will give a disjoint union of irreducible components, each isomorphic to $\mathbb{P}^{2}$, in the special fiber of the special level Shimura variety of type $(2,1)$.
  This special level Shimura variety has been studied in \cite[\S 5]{badreduction}.
\end{remark}
\subsubsection{The second isogeny}
Next, we consider the quotient $\omega_{\mathcal{A}_1} \twoheadrightarrow \mathcal{F}_2 / \mathcal{F}_1$ (the map is the obvious one from \eqref{eq:decomp-first-vb}).
For the same reason as before, the quotient is $V_{\mathcal{A}_1}$-stable, and we obtain
an $F, V$-kernel $G_1 \hookrightarrow \mathcal{A}_1$. Let $\mathcal{A}_2$ denote $\mathcal{A}_1 / G_1$.
Copying the argument of the previous section gives us the following diagram
with exact rows
\begin{equation} \label{eq:hodge-filtr-2}
\begin{tikzcd}
0 & {\mathcal{G}^{2, \bot}/\mathcal{G}^{1, \bot}} & {\omega_{\mathcal{A}_2}} & {\omega_{\mathcal{A}_1}} & {\mathcal{F}_2/\mathcal{F}_1} & 0 \\
0 & {\mathcal{G}^{2, \bot}/\mathcal{G}^{1, \bot}} & {\Hdr(\mathcal{A}_2)} & {\Hdr(\mathcal{A}_1)} & {\mathcal{G}^{2, \bot}/\mathcal{G}^{1, \bot}} & 0 \\
  0 & {\mathcal{F}_2/\mathcal{F}_1} & {\Lie(\mathcal{A}_2^\vee)} & {\Lie(\mathcal{A}_1^\vee)} & {\mathcal{G}^{2, \bot}/\mathcal{G}^{1, \bot}} & 0
	\arrow[from=1-1, to=1-2]
	\arrow[from=1-2, to=1-3]
	\arrow["=", from=1-2, to=2-2]
	\arrow[from=1-3, to=1-4]
	\arrow[from=1-3, to=2-3]
	\arrow[from=1-4, to=1-5]
	\arrow[from=1-4, to=2-4]
	\arrow[from=1-5, to=1-6]
	\arrow["0", from=1-5, to=2-5]
	\arrow[from=2-1, to=2-2]
	\arrow[from=2-2, to=2-3]
	\arrow["0", from=2-2, to=3-2]
	\arrow[from=2-3, to=2-4]
	\arrow[from=2-3, to=3-3]
	\arrow[from=2-4, to=2-5]
	\arrow[from=2-4, to=3-4]
	\arrow[from=2-5, to=2-6]
	\arrow["=", from=2-5, to=3-5]
	\arrow[from=3-1, to=3-2]
	\arrow[from=3-2, to=3-3]
	\arrow[from=3-3, to=3-4]
	\arrow[from=3-4, to=3-5]
	\arrow[from=3-5, to=3-6]
\end{tikzcd}
\end{equation}
We can also read off the following decompositions into isotypic components
\begin{equation} \label{eq:decomp-second-vb-1}
  \omega_{\mathcal{A}_2} \cong \omega_{\mathcal{A}_2, \tau_0} \oplus \omega_{\mathcal{A}_2, \tau_1} \cong \mathcal{F}_1 \oplus \mathcal{G}^{3, \bot} / \mathcal{G}^{1, \bot}
\end{equation}
and
\begin{equation} \label{eq:decomp-second-vb-2}
  \Lie(\mathcal{A}_2^\vee) \cong \Lie(\mathcal{A}_2^\vee)_{\tau_0} \oplus \Lie(\mathcal{A}_2^\vee)_{\tau_1} \cong \mathcal{F}_3 / \mathcal{F}_1 \oplus \mathcal{G}^{1, \bot}.
\end{equation}
In particular, we see that $\mathcal{A}_2$ has type $(1, 2)$.
Finally, by viewing $\mathcal{F}_3 / \mathcal{F}_1$ as a quotient of $\Hdr(\mathcal{A}_2)$ via
\eqref{eq:decomp-second-vb-2}, we get the following induced isomorphisms by combining Lemmas \ref{lem:frob-ver-descr} and \ref{lem:frob-ver-reln}.
\[
  F_{\mathcal{A}_2} : (\mathcal{F}_3 / \mathcal{F}_1)^{(p)} \xlongrightarrow{\sim}
  \mathcal{G}^{3, \bot}/\mathcal{G}^{1, \bot} \subseteq \Hdr(\mathcal{A}_2)
\]
and
\begin{equation} \label{eq:VX2-descr}
  V_{\mathcal{A}_2} : \mathcal{F}_3 / \mathcal{F}_1 \xlongrightarrow{\sim}
  (\mathcal{G}^{3, \bot}/\mathcal{G}^{1, \bot})^{(p)} \subseteq \Hdr(\mathcal{A}_2)^{(p)}.
\end{equation}

\subsubsection{The third isogeny}
Let us now restrict to the closed stratum $\Zcb'_{21}$. Its defining condition
\eqref{eq:z21'-cond} implies that we have $\mathcal{G}_1^\bot
\subseteq \mathcal{F}_1^\bot$. Thus, by \eqref{eq:decomp-second-vb-1}, we can consider the following quotient
\begin{equation} \label{eq:quot-isog-3}
  \omega_{\mathcal{A}_2|_{\Zcb'_{21}}} \twoheadrightarrow \mathcal{G}^{3, \bot} / \mathcal{F}_1^{\bot}.
\end{equation}
\begin{lemma}
  The kernel of \eqref{eq:quot-isog-3} is $V_{\mathcal{A}_2}$-stable.
\end{lemma}
\begin{proof}
  Since $V_{\mathcal{A}_2}$ exchanges the isotypic
  components, it suffices to show that 
  $\omega_{\mathcal{A}_2, \tau_0} \subseteq \ker V_{ \mathcal{A}_2}$.
  Since $\mathcal{A}_2$ is of type $(1, 2)$, we have that $(\ker V_{ \mathcal{A}_2})_{\tau_0}$ must be of rank $1$.
  It then follows from \eqref{eq:VX2-descr} that $(\ker V_{ \mathcal{A}_2})_{\tau_0} = \omega_{
\mathcal{A}_2, \tau_0}$, so we are done.
\end{proof}
We can now apply Lemma \ref{lem:integrate-subgroup}
to obtain an $\mathcal{O}_{E}$-stable subgroup $G_3$ of $\mathcal{A}_2|_{\Zcb'_{21}}$ 
of order $p$ corresponding
to \eqref{eq:quot-isog-3}.
Let $\mathcal{A}_2|_{\Zcb'_{21}} \rightarrow \mathcal{A}_3$, denote the corresponding isogeny.
Moreover, Lemma \ref{lem:ht-one} implies that $G$ must be an $F, V$-kernel. Thus,
we can apply Lemma \ref{lem:hodge-seq-final} to determine the Hodge filtration as
follows
\begin{equation} \label{eq:hodge-filtr-3}
  \begin{tikzcd}
  0 & {(\mathcal{G}^{3, \bot} / \mathcal{F}_1^{\bot})^{(p)}} & {\omega_{\mathcal{A}_3}} & {\omega_{\mathcal{A}_2|_{\Zcb'_{21}}}} & {\mathcal{G}^{3, \bot} / \mathcal{F}_1^{\bot}} & 0 \\
    0 & {\mathcal{H}^0(\Delta(G_3)_{\Zcb'_{21}})} & {\Hdr(\mathcal{A}_3)} & {\Hdr(\mathcal{A}_2|_{\Zcb'_{21}})} & {\mathcal{D}(G_3)_{\Zcb'_{21}}} & 0 \\
	0 & {\mathcal{G}^{3, \bot} / \mathcal{F}_1^{\bot}} & {\Lie(\mathcal{A}_3^\vee)} & {\Lie(\mathcal{A}_2^\vee|_{\Zcb'_{21}})} & {(\mathcal{G}^{3, \bot} / \mathcal{F}_1^{\bot})^{(p)}} & 0
	\arrow[from=1-1, to=1-2]
	\arrow[from=1-2, to=1-3]
	\arrow["\simeq", from=1-2, to=2-2]
	\arrow[from=1-3, to=1-4]
	\arrow[from=1-3, to=2-3]
	\arrow[from=1-4, to=1-5]
	\arrow[from=1-4, to=2-4]
	\arrow[from=1-5, to=1-6]
	\arrow["0", from=1-5, to=2-5]
	\arrow[from=2-1, to=2-2]
	\arrow[from=2-2, to=2-3]
	\arrow["0", from=2-2, to=3-2]
	\arrow[from=2-3, to=2-4]
	\arrow[from=2-3, to=3-3]
	\arrow[from=2-4, to=2-5]
	\arrow[from=2-4, to=3-4]
	\arrow[from=2-5, to=2-6]
	\arrow["\simeq", from=2-5, to=3-5]
	\arrow[from=3-1, to=3-2]
	\arrow[from=3-2, to=3-3]
	\arrow[from=3-3, to=3-4]
	\arrow[from=3-4, to=3-5]
	\arrow[from=3-5, to=3-6]
\end{tikzcd}
\end{equation}
The following lemma gives a better description of $\mathcal{D}(G_3)_{\Zcb'_{21}}$ as a quotient
of $\Hdr(\mathcal{A}_2|_{\Zcb'_{21}})$.
\begin{lemma} \label{lem:ver-a2}
  View $\mathcal{F}_3 / \mathcal{G}^1$ as a quotient of $\Hdr( \mathcal{A}_2|_{\Zcb'_{21}})$ via \eqref{eq:decomp-second-vb-2}.
  Then $V_{ \mathcal{A}_2} : \Hdr(\mathcal{A}_2|_{\Zcb'_{21}}) \twoheadrightarrow
  \omega_{\mathcal{A}_2|_{\Zcb'_{21}}}^{(p)}$ induces an isomorphism between quotients
  $\mathcal{F}_3 / \mathcal{G}^1 \xrightarrow{\sim} (\mathcal{G}^{3, \bot} /
  \mathcal{F}_1^{\bot})^{(p)}$.
  In particular we can identify $\mathcal{F}_3 / \mathcal{G}^1$ with $\mathcal{D}(G_3)_{\Zcb'_{21}}$, as quotients of of $\Hdr(\mathcal{A}_2|_{\Zcb'_{21}})$.
\end{lemma}
\begin{proof}
  The second part follows from the first since, by Lemma \ref{lem:hodge-seq-final}, the isomorphism $\mathcal{D}(G_3)_{\Zcb'_{21}} \cong (\mathcal{G}^{3, \bot} / \mathcal{F}_1^{\bot})^{(p)}$ in the above diagram
  is also induced by $V_{\mathcal{A}_2}$.
  We will show the first part by applying Lemma \ref{lem:frob-ver-descr}
  to a \textit{different}
sequence of isogenies $\mathcal{A}' \rightarrow
\mathcal{A}_1' \rightarrow \mathcal{A}_2$ (this sequence will only be defined over ${\Zcb'_{21}}$). Let $G_0'$ denote the kernel of $\mathcal{A}' \rightarrow \mathcal{A}_1'$, and $G_1'$ denote the kernel of $\mathcal{A}_1' \rightarrow \mathcal{A}_2$.
The point is that we will pick $G_0'$ such that
$\mathcal{D}(G_0')_{\Zcb'_{21}} \cong \mathcal{G}^{3, \bot} / \mathcal{F}_1^{\bot}$ which will allow
us to use
Lemma \ref{lem:frob-ver-descr} to understand its image under $V_{ \mathcal{A}_1}$ (and
then also under $V_{ \mathcal{A}_2}$).
Let us carry out this strategy now. We will be brief since the construction of these
isogenies is
essentially identical to what has been done above (i.e. apply Lemmas \ref{lem:integrate-subgroup} and \ref{lem:hodge-seq-final}).
All the following will take place over ${\Zcb'_{21}}$.
Let $\mathcal{A}' \rightarrow \mathcal{A}_1'$ be obtained from the quotient
$\omega_{\mathcal{A}'} \rightarrow \mathcal{F}_3 /
\mathcal{G}^1$, and then note that $\mathcal{A}_1' \rightarrow \mathcal{A}_2$ can be
obtained from the quotient $\omega_{\mathcal{A}_1'} \rightarrow \mathcal{G}^1 /
\mathcal{F}_1$ (indeed, the composition of the isogenies corresponds to the quotient $\omega_{\mathcal{A}'}
\rightarrow \mathcal{F}_3 / \mathcal{F}_1$). By Lemma \ref{lem:frob-ver-reln} iii) we
have that $F_{\mathcal{A}'}(\mathcal{G}^{3, \bot, (p)}) / F_{\mathcal{A}'}(\mathcal{F}_1^{\bot, (p)}) = \mathcal{F}_3 /
\mathcal{G}^1$. Treating $\mathcal{F}_3 /
\mathcal{G}^1$ as a subquotient of $\Hdr(\mathcal{A}'_1)$ and applying
Lemma \ref{lem:frob-ver-descr} then shows that
$V_{\mathcal{A}_1'}$ induces an isomorphism
\[
  \mathcal{F}_3 / \mathcal{G}^1 \xrightarrow{\sim} \mathcal{H}^0(\Delta(G_1')_{\Zcb'_{21}})^{(p)} \cong (\mathcal{G}^{3, \bot} / \mathcal{F}_1^{\bot})^{(p)}.
\]
Finally, one observes that both $\mathcal{F}_3 / \mathcal{G}^1$ and $\mathcal{G}^{3, \bot} /
\mathcal{F}_1^{\bot}$ can be viewed as subquotients of the image of
$\Hdr(\mathcal{A}_2|_{\Zcb'_{21}})$ in $\Hdr(\mathcal{A}'_1)$, and hence we get that $V_{\mathcal{A}_2}$ also induces such an isomorphism.
\end{proof}
Let $\iota_3 : \mathcal{O}_E \rightarrow \End(\mathcal{A}_3)$ denote the induced $\mathcal{O}_E$-action
on $\mathcal{A}_3$. We can now prove the following theorem.
\begin{theorem} \label{thm:p-div-family}
    The abelian variety $\mathcal{A}_3/{\Zcb'_{21}}$ is of type $(2, 1)$ with respect to $\iota_3$, and it also inherits a prime-to-$p$ quasi-polarization
    $\lambda_3$ from $\mathcal{A}'$.
    Furthermore, $\omega_{\mathcal{A}_{3, \tau_0}}$ obtains a natural filtration
    \[
      0 \longrightarrow \mathcal{F}_3 / \mathcal{G}^1 \longrightarrow \omega_{\mathcal{A}_3, \tau_0}
        \longrightarrow \mathcal{F}_1 \longrightarrow 0
    \]
    and we also get
    \[
      \omega_{\mathcal{A}_3, \tau_1} \cong \mathcal{F}_1^{\bot} / \mathcal{G}^{1, \bot}.
  \]
    The isotypic components of $\Lie(\mathcal{A}_3^\vee)$ can similarly be described via
    \[
      \Lie(\mathcal{A}_3^\vee)_{\tau_0} \cong \mathcal{G}^1 / \mathcal{F}_1
    \]
    and
    \[
    0 \longrightarrow \mathcal{G}^{3, \bot} / \mathcal{F}_1^{\bot} \longrightarrow \Lie(\mathcal{A}_3^\vee)_{\tau_1} \longrightarrow \mathcal{G}^{1, \bot} \longrightarrow 0.
    \]
\end{theorem}
\begin{proof}
  The description of isotypic components easily follows from \eqref{eq:decomp-second-vb-1}, \eqref{eq:decomp-second-vb-2} and \eqref{eq:hodge-filtr-3}.
    This also implies that the signature of $\mathcal{A}_3$ is $(2, 1)$.
    Let $\lambda'$ denote the prime-to-$p$ quasi-polarization associated to $\mathcal{A}'$.
    We will show that $p\lambda'$ on $\mathcal{A}'$ descends to a
    prime-to-$p$ quasi-polarization on $\mathcal{A}_3$.
    Let us use the notation $\Phi : \mathcal{A}'|_{\Zcb'_{21}} \rightarrow \mathcal{A}_3$.
    Consider the short exact sequence
    \[
      0 \longrightarrow \Phi(\Hdr(\mathcal{A}_3)) \longrightarrow \Hdr(\mathcal{A}'|_{\Zcb'_{21}}) \longrightarrow \mathcal{D}(\ker \Phi)_{\Zcb'_{21}} \longrightarrow 0.
    \]
    If we can show that $\Phi(\Hdr(\mathcal{A}_3))$
    is maximally isotropic with respect to the non-degenerate pairing $\lambda'$ on
    $\Hdr(\mathcal{A}'|_{\Zcb'_{21}})$, we would be done,
    as then $p\lambda'$ must descend to a prime-to-$p$ quasi-polarization on
    $\mathcal{A}_3$.
    Now, the morphism $\Hdr(\mathcal{A}_3)
    \rightarrow \Hdr(\mathcal{A}'|_{\Zcb'_{21}})$ is nothing
    but the following composition of morphisms
    \[
      \Hdr(\mathcal{A}_3)
        \longrightarrow \Hdr(\mathcal{A}_2|_{\Zcb'_{21}})
        \longrightarrow \Hdr(\mathcal{A}_1|_{\Zcb'_{21}})
        \longrightarrow \Hdr(\mathcal{A}'|_{\Zcb'_{21}}).
    \]
    Hence, tracing through diagrams \eqref{eq:hodge-filtr-1}, \eqref{eq:hodge-filtr-2} and \eqref{eq:hodge-filtr-3} (and applying Lemma \ref{lem:ver-a2}), we see that the
    image can be identified with $\mathcal{G}^1 \oplus \mathcal{G}^{1, \bot}$.
    This is clearly maximally isotropic for $\lambda'$, so we are done.
\end{proof}
\begin{remark} \label{rmk:iwahori}
    As mentioned earlier, one can show that the morphism $\mathcal{A}_1 \rightarrow \mathcal{A}_3$
    gives us a $\Zcb'_{21}$-valued point of the corresponding Iwahori level Shimura variety. It is possible to use this to identify ${\Zcb'_{21}}$ with a disjoint union of  irreducible components at Iwahori
    level. This gives a description of the ``supersingular screens'' of \cite[\S 4]{badreduction} in terms of Deligne--Lusztig varieties.
\end{remark}
\begin{proposition} \label{prop:map-to-fl}
The abelian variety $\mathcal{A}_3 / \Zcb'_{21}$ induces a morphism $\Zcb'_{21} \rightarrow
\FlKtwo$.
\end{proposition}
\begin{proof}
  From Theorem \ref{thm:p-div-family}, we are already have a triple $\mathcal{(A}_3, \iota_3, \lambda_3)$.
  Furthermore, the theorem also gives us a refinement of the Hodge filtration given by
  \begin{equation} \label{eq:refined-hodge-filtr}
    0 \subseteq \mathcal{F}_3 / \mathcal{G}^1 \subseteq \omega_{ \mathcal{A}_3, \tau_0} \subseteq
  \Hdr( \mathcal{A}_3)_{\tau_0}.
  \end{equation}
  Thus, to obtain a morphism $\Zcb'_{21}\rightarrow \FlKtwo$, it remains to equip $\mathcal{A}_3$ with
  a level structure. Since the isogeny is of $p$-power degree, 
  we can simply transfer the $K'^p$-level structure of $\mathcal{A}'$
  to a $K^p$-level structure of $\mathcal{A}_3$ via the isogeny we just constructed.
\end{proof}
\subsection{Identification of the closed strata}
The morphism of Proposition \ref{prop:map-to-fl} will not be injective, but its restriction to $\Zcb'_1$ will be.
Therefore, let $i : \Zcb'_1 \rightarrow \FlKtwo$ denote the restriction of this morphism to $\Zcb'_1$.
The goal of this subsection is to use $i$ to identify $\Zcb'_1$ with $\Zcb_{2}$.

Note that \eqref{eq:refined-hodge-filtr} allows us to treat $\mathcal{F}_3 /
\mathcal{G}^1$ as a subbundle of $\Hdr(\mathcal{A}_3)_{\tau_0}$. We can then consider the orthogonal complement, with respect to $\lambda_3$, $(\mathcal{F}_3 /
\mathcal{G}^1)^\bot$ as a subbundle of $\Hdr(\mathcal{A}_3)_{\tau_1}$.
\begin{lemma}
  As subbundles of $\Hdr(\mathcal{A}_3)_{\tau_0}$, we have that $F_{\mathcal{A}_3}((\mathcal{F}_3 / \mathcal{G}^1)^{\bot, (p)}) = \mathcal{F}_3 / \mathcal{G}^1$.
  In particular, $i$ factors through $\Zcb_{2}$.
\end{lemma}
\begin{proof}
  The second part follows from the first by using the description of the Hodge filtration of $\Hdr(
  \mathcal{A}_3)_{\tau_0}$ from \eqref{eq:refined-hodge-filtr}, and the description
  of the refined conjugate filtration from \eqref{eq:conj-filtr-21-def}.

  Let $\mathcal{E}$ denote the pre-image of $\mathcal{G}^{3, \bot} / \mathcal{F}_1^{\bot}$ under the morphism $\Hdr(\mathcal{A}_3)_{\tau_1} \to \Lie(\mathcal{A}_3^\vee)_{\tau_1}$.
  This is an extension of $\mathcal{G}^{3, \bot} / \mathcal{F}_1^{\bot}$ by
  $\omega_{\mathcal{A}_3, \tau_1}$.
  It follows from Lemmas \ref{lem:ver-a2} and \ref{lem:frob-ver-descr}
that $F_{\mathcal{A}_3}(\mathcal{E}^{(p)}) = \mathcal{F}_3 /
  \mathcal{G}^1$.
  Tracing through the diagrams \eqref{eq:hodge-filtr-1}, \eqref{eq:hodge-filtr-2} and \eqref{eq:hodge-filtr-3},
  one sees that $\mathcal{F}_3 / \mathcal{G}^1 \oplus \mathcal{E}$ is exactly the kernel
  of $\Hdr(\mathcal{A}_3) \rightarrow \Hdr(\mathcal{A}'|_{\Zcb'_1})$,
  and since the image is totally isotropic with respect to $\lambda'$, we must have that
  $(\mathcal{F}_3 / \mathcal{G}^1)^\bot = \mathcal{E}$ in $\Hdr(\mathcal{A}_3)$.
\end{proof}
\begin{proposition} \label{prop:strata-iso}
  The map $i : \Zcb'_{1} \rightarrow \Zcb_{2}$ is an isomorphism.
\end{proposition}
\begin{proof}
  It is possible to construct an inverse to $i$, by performing the same process
  as in §\ref{sect:exotic-hecke}, but starting with the universal family $\mathcal{A} / \Zcb_{2}$.
  Since the fact that this is an isomorphism should be no surprise, given the results in \cite{vollaard},
  let us only sketch this process, and leave the details to the reader.
  First, on the stratum $\Zcb_{21}$, we can consider the rank $2$ quotient $\omega_{\mathcal{A}_3} \rightarrow \omega_{\mathcal{A}_3} /
  \mathcal{F}_1$.
  This gives us a degree $p^2$-isogeny $\mathcal{A} \rightarrow \mathcal{A}_1'$.
  It is clear that $\mathcal{A}_1'$ is still of type $(2, 1)$.
  The kernel of $\mathcal{A} \rightarrow \mathcal{A}_1'$ is not an $F, V$-kernel (not even on $\Zcb_{2}$),
  but one can still check pointwise, using Dieudonné theory, that the quotient
  $\omega_{\mathcal{A}_1'} \rightarrow \omega_{\mathcal{A}_1', \tau_1}$ is
  $V$-stable over $\Zcb_{2}$.
  Thus, this quotient defines an isogeny $\mathcal{A}_1' \rightarrow \mathcal{A}_2'$.
  It is clear that $\mathcal{A}_2'$ will have type $(3, 0)$, and one checks that
  it furthermore obtains a prime-to-$p$ quasi-polarization (by showing that the image of $\Hdr(\mathcal{A}_2')$ in $\Hdr(\mathcal{A})$ is $\mathcal{F}_1 \oplus V^{-1}(\mathcal{F}_1^{(p)}) = \mathcal{F}_1 \oplus \mathcal{F}_1^\bot$), and that the induced filtration satisfies the defining
  condition of $\Zcb'_1$ (this is easiest to do pointwise using Dieudonné modules). This gives us a map $\Zcb_{2} \rightarrow \Zcb'_{1}$.
  
  Next, one notices that by applying the construction of §\ref{sect:exotic-hecke} to $\mathcal{A}_2'$
  and composing the two isogenies, one gets the multiplication by $p$ isogeny. The same works in the other direction.
  This implies that the map $\Zcb_{2} \rightarrow \Zcb'_{1}$ is the inverse of $\Zcb'_{1} \rightarrow
  \Zcb_{2}$, so we are done.
  \end{proof}

\subsection{The automorphic line bundles}
From Theorem \ref{thm:p-div-family}, we have the following isomorphism (where we
identified $\Zcb'_1 \cong \Zcb_2$ via $i$)
\begin{equation} \label{eq:L-descr}
  \mathcal{L}(a, b, c, d)|_{\Zcb'_1} \cong \mathcal{F}_1^{\otimes a} \otimes (\mathcal{F}_3 /
  \mathcal{G}^1)^{\otimes b} \otimes (\mathcal{G}^1 / \mathcal{F}_1)^{\otimes c} \otimes \delta^d.
\end{equation}
By using the equality $\mathcal{G}^1 = \mathcal{F}_2$ on $\Zcb'_1$, we can easily compute this
as
\[
  \mathcal{L}(a, b, c, d)|_{\Zcb'_1} \cong \mathcal{L}'(b, c, a, d)|_{\Zcb'_1}.
\]
A problem with this description is that if $(a, b, c, d)$ is dominant, then
$(b, c, a, d)$ will usually not be.
Furthermore, as we mentioned previously, it will be easier compute the global sections if
we can view $\mathcal{L}(b, c, a, d)|_{\Zcb'_1}$ as the restriction of some
line bundle on $\Fl'_{K', P^*}$.
Recall from \eqref{eq:invertible-hasse} that we have $\mathcal{O}_{\Zcb'_1} \cong
\mathcal{L}'(1, 0, p, -p)$. Hence, shifting by $-(b - c)$ times this weight gives us
\[
  \mathcal{L}'(b, c, a, d)|_{\Zcb'_1} \cong \mathcal{L}'(c, c, a - p(b - c), (b - c)p + d)|_{\Zcb'_1} 
\]
This second line bundle can be obtained via restriction from $\Fl_{K', P^*}$.
We have proved the following theorem.
\begin{theorem} \label{thm:aut-vb-descr}
    We have an isomorphism
    \begin{align*}
      \mathcal{L}(a, b, c, d)|_{\Zcb'_1} \cong \mathcal{L}'_{P^*}(c, c, a - p(b - c), (b -
      c)p + d)|_{\Zcb'_1}.
    \end{align*}
    In particular, for the irregular weights $(a, b, b - 2, d)$ we get
    \[
      \mathcal{L}(a, b, b - 2, d)|_{\Zcb'_1} \cong \mathcal{L}'_{P^*}(b - 2, b - 2, a - 2p, 2p + d)|_{\Zcb'_1}.
    \]
\end{theorem}

\section{The global sections} \label{sect:global-sections}
Let $p : \Zcb'_1 \rightarrow \ShKtr$ denote the composition $\Zcb'_1 \rightarrow \Fl'_{K', P^*} \rightarrow
\ShKtr$. We will also treat $\Zcb'_1$ as a closed subscheme of $\ShKtwo$ via $i$.
Our goal in this section is to prove the following theorem, which can be seen as describing the global sections of $\mathcal{L}(\lambda)|_{\Zcb'_1}$, for a $p$-restricted irregular weight $\lambda$,
in terms of spaces of algebraic modular forms with explicit regular weights. Here, for a
character $\lambda \in X^*(T)$, we let $W(\lambda)$ denote the corresponding dual Weyl module of $G_{\Fp}$,
and we let $L(\lambda)$ denote the irreducible representation of highest weight $\lambda$ (which is
the socle of $W(\lambda)$).
\begin{theorem} \label{thm:main-decomp}
  Given an irregular weight $\lambda = (a, b, b - 2, d) \in X^*(T)$, set $\lambda_1 := (p + b - 1,
  a, b, d + 2)$ and $\lambda_2 = (p + b, p + b - 1, a + 1, d + 1)$. Then, if $a - b \le p - 2$,
  we have the following exact sequence
  \[
      0 \longrightarrow \omega'(W(\lambda_2)) \longrightarrow p_*(\mathcal{L}(\lambda)|_{\Zcb'_1}) \longrightarrow \omega'(L(\lambda_1)) \longrightarrow 0.
  \]
  For $a - b = p - 1$ we get
  $p_*(\mathcal{L}(\lambda)|_{\Zcb'_1}) = \omega'(L(\lambda_1))$.
\end{theorem}
Let us first explain how this reduces to a proposition about the cohomology
of line bundles on a Deligne--Lusztig variety.
Recall from Remark \ref{rmk:DL-var} that we have the $1$-dimensional closed Deligne--Lusztig variety $\overline{X}_1 \subseteq \mathbb{P}^{2, *}$.
It fits into the following diagram, where all squares are cartesian
\[\begin{tikzcd}
	{\Zcb'_1} & {\overline{\mathcal{C}}'_1} & {\overline{X}_1} \\
  {\Fl'_{K', P*}} & {\Gzf^{\mu', P^*}} & {G / P^* \cong \mathbb{P}^{2, *}} \\
  \ShKtr & {\Gzip^{\mu'}} & {*.}
	\arrow[from=1-1, to=1-2]
	\arrow[from=1-1, to=2-1]
	\arrow[from=1-2, to=2-2]
	\arrow[from=1-3, to=1-2]
	\arrow[from=1-3, to=2-3]
	\arrow[from=2-1, to=2-2]
	\arrow[from=2-1, to=3-1]
	\arrow[from=2-2, to=3-2]
	\arrow[from=2-3, to=2-2]
	\arrow[from=2-3, to=3-3]
	\arrow[from=3-1, to=3-2]
	\arrow[from=3-3, to=3-2]
\end{tikzcd}\]
In particular, by pulling back from $\Gzf^{\mu', P^*}$ we get a line bundle
$\mathcal{O}_{\mathbb{P}^{2, *}}(\lambda)$ on $\mathbb{P}^{2, *}$ for each character $\lambda = (a, a, b, d) \in X^*(T)$. We will let $\mathcal{O}_{\overline{X}_1}(\lambda)$ denote the restriction $\mathcal{O}_{\mathbb{P}^{2, *}}(\lambda)|_{\overline{X}_1}$.
Furthermore, it follows from Lang's theorem \cite[Theorem 10.1]{steinberg} that $\Gzip^{\mu'} \cong [* / \underline{G(\Fp)}]$.
Hence, the horizontal arrows on the right are all $\underline{G(\Fp)}$-torsors. In particular
$\mathcal{O}_{\overline{X}_1}(\lambda)$ is $G(\Fp)$-equivariant. We are now ready to state the
following.
\begin{proposition} \label{prop:main-decomp}
  In the notation of the previous theorem, set $\lambda' := (b - 2, b - 2, a - 2p, 2p + d)$ (this is the weight corresponding to $\lambda$ via Theorem \ref{thm:aut-vb-descr}). Then, if $a - b \le p - 2$, we have the following exact sequence of $G(\Fp)$-representations
  \[
    0 \longrightarrow W(\lambda_2) \longrightarrow H^0(\overline{X}_1, \mathcal{O}_{\overline{X}_1}(\lambda')) \longrightarrow L(\lambda_1) \longrightarrow 0.
  \]
  For $a - b = p - 1$ we get
  $H^0(\overline{X}_1, \mathcal{O}_{\overline{X}_1}(\lambda)) = L(\lambda_1)$, as
  $G(\Fp)$-representations.
\end{proposition}
\begin{proof}[Proof that Proposition \ref{prop:main-decomp} implies Theorem \ref{thm:main-decomp}]
Since the automorphic vector bundles can be obtained via pullback from $\Gzip^{\mu'}$,
it follows from flat base change for algebraic stacks \cite[Proposition 7.1.14]{alper} (using that $\ShKtr \to \Gzip^{\mu'}$ is flat) that it suffices to compute the corresponding pushforward along the middle column in the above diagram.
Furthermore, by descent, it suffices to compute the pushforward on the right-hand column
as $G(\Fp)$-equivariant vector bundles.
Now, Theorem \ref{thm:aut-vb-descr} implies that 
the line bundle on $\overline{X}_1$ corresponding to $\mathcal{L}(\lambda)|_{\Zcb'_1}$ is $\mathcal{O}_{\overline{X}_1}(\lambda')$, and $W(\lambda_2)$ resp. $L(\lambda_1)$ clearly corresponds to $\omega'(W(\lambda_2))$ resp. $\omega'(L(\lambda_1))$ on $\ShKtr$. So we are done.
\end{proof}
It therefore remains to prove Proposition \ref{prop:main-decomp}, which we will do in \S \ref{sect:global-sect-zcb}.

Let us now make a temporary change of setting.
Namely, it is somewhat unfortunate that the partial flag variety corresponding to 
$\Fl_{K, P^*}$ is $\mathbb{P}^{2, *}$ and not $\mathbb{P}^2$ (this will introduce a bunch of minus signs in our computations, as the twisting sheaf $\mathcal{O}_{\mathbb{P}^{2, *}}(n)$ has weight $(0, 0, -n, 0)$ as opposed to $(n, 0, 0, 0)$).
Hence, to avoid excessive minus signs, let us instead work 
on $\mathbb{P}^2$.
We will denote the corresponding closed Deligne--Lusztig variety by $\overline{X}_2$ (which
corresponds to $s_2 \in W(T)$).
This is the same as the Deligne--Lusztig variety considered in \cite{vollaard}, and it follows
from \cite[Remark 4.11b)]{vollaard} that $\overline{X}_2$ can be identified with the Fermat curve cut out
by the equation $x^{p + 1} + y^{p + 1} + z^{p + 1}$.
The next two subsections will be devoted to computing global sections of line bundles on
$\overline{X}_2$. We will then transfer the results back to $\overline{X}_1$, via the isomorphism $G / P
\cong G / P^*$ which sends $g \mapsto w_0g^{-t}w_0$. This has the effect
$(a, b, c, d) \mapsto (-c, -b, -a, -d)$ on weights.

As above, we will use the notation $\mathcal{O}_{\mathbb{P}^2}(\lambda)$ and
$\mathcal{O}_{\overline{X}_2}(\lambda)$ for the line bundle corresponding to the character $\lambda
= (a, b, b, d) \in X^*(T)$.
Note that $\mathcal{O}_{\mathbb{P}^2}(k, 0, 0, 0) \cong \mathcal{O}_{\mathbb{P}^2}(k)$ (the usual twisting sheaf).
In the next two subsections, we will for simplicity only work with $\mathcal{O}_{\mathbb{P}^2}(k)$,
as the general case will follow from twisting by the determinant and multiplier characters.
We will also use the notation $\mathcal{O}_{\overline{X}_2}(k) := \mathcal{O}_{\overline{X}_2}(k, 0,
0, 0)$.

\subsection{Global sections over $\mathbb{P}^2$}
Let us first recall some results of \cite{doty}, describing the
structure of $H^0(\mathbb{P}^n, \mathcal{O}_{\mathbb{P}^n}(k))$ as an algebraic $\GL_n$-representation over $\Fpb$.
Here, $\mathbb{P}^n$ denotes the projective space parametrizing $1$-dimensional quotients, and $\mathcal{O}_{\mathbb{P}^n}(k)$ is the usual twisting sheaf.
We will immediately specialize these results to the case $n = 2$ and $k = ap + b$,
for $0 \le a, b < p$.

In this subsection we will be working over $\Fpb$.

\subsubsection{$p$-carry patterns}
\begin{definition}
  Fix $k = ap + b$, with $0 \le a, b < p$. Given a monomial $x^{\lambda_0}y^{\lambda_1}z^{\lambda_2} \in H^0(\mathbb{P}^2,
  \mathcal{O}_{\mathbb{P}^2}(k))$, let $\lambda_{i, 0}$ denote the $p^0$-digit in the $p$-digit
  expansion of $\lambda_i$ (for $i \in \{0, 1, 2\}$).
  We say that the monomial has ``$p$-carry pattern'' $c(\underline{\lambda}) \in \mathbb{Z}_{\ge 0}$ (or just denoted $c$ when no confusion can occur), if
  $\lambda_{0, 0} + \lambda_{1, 0} + \lambda_{2, 0} = c(\underline{\lambda})p + r$ for some remainder $0 \le r < p$.
  Let $C(k) \subseteq \mathbb{Z}_{\ge 0}$ denote the set of carry patterns that can arise from a degree $k$ monomial.
  It inherits a natural ordering coming from the integers.
\end{definition}
Note that since we are working with monomials of $3$ variables, we have $C(k) \subseteq
      \{0, 1, 2\}$ for all $k$ as in the definition.
\begin{remark}
    In \cite{doty}, the above is defined for any $n$ and any $k \ge 0$. Instead of an integer $c$, one
    obtains a tuple measuring the carrying in each $p^i$-digit in the $p$-digit expansion of $k$. We
    have specialized it above in order to simplify the notation in what follows.
\end{remark}
\subsubsection{Global sections}
The following theorem is a specialization of \cite[Theorem]{doty}
\begin{theorem} \label{thm:subreps-H0}
    Fix $k = ap + b$, with $0 \le a, b < p$. We have an ordered bijection
    \[
      \{ \text{sub-representations of } H^0(\mathbb{P}^2, \mathcal{O}_{\mathbb{P}^2}(k)) \} \longleftrightarrow C(k).
    \]
    It is explicitly given by sending $c \in C(k)$ to the subrepresentation of monomials of
    carry-pattern $\le c$.
    In particular, the lattice of subrepresentations is totally ordered, and hence
    $H^0(\mathbb{P}^2, \mathcal{O}_{\mathbb{P}^2}(k))$ comes with a natural increasing filtration, which we
    will call the \textit{carry-filtration}.
\end{theorem}
By an easy calculation one deduces the following lemma,
which describes the length of the carry-filtration
of $H^0(\mathbb{P}^2, \mathcal{O}_{\mathbb{P}^2}(k))$.
\begin{lemma}
    Fix $k = ap + b$, with $0 \le a, b < p$. Then we have
    \[
      C(k) =
      \begin{cases}
        \{ 0 \} & \text{if } a = 0,                    \\
        \{ 0, 1 \} & \text{if } a = 1 \text{ or } (2 \le a \text{ and } p - 3 < b), \\
        \{ 0, 1, 2 \} & \text{if } 2 \le a \text{ and } b \le p - 3.
      \end{cases}
    \]
\end{lemma}
Let us now describe the irreducible constituents of the carry-filtration explicitly.
Note that, we can identify $H^0(\mathbb{P}^2, \mathcal{O}_{\mathbb{P}^2}(a))$ with the dual Weyl module $W(a, 0, 0, 0)$. 
\begin{corollary} \label{cor:weyl-decomp}
  Fix $k = ap + b$, with $0 \le a, b < p$. Then the carry-filtration of the $G_{\Fpb}$-representation $W(ap + b, 0, 0, 0)$ has the following, ordered,
list of irreducible constituents
  \begin{align*}
    & \begin{cases}
      \left\{ L(ap + b, 0, 0, 0) \right\} & \text{if } a = 0,                    \\
      \left\{ L(ap + b, 0, 0, 0), L((a - 1)p + p - 1, b + 1, 0, 0) \right\} & \text{if } a = 1 \text{ or } (2 \le a \text{ and } p - 3 < b), \\
			\{ L(ap + b, 0, 0, 0), L((a - 1)p + p - 1, b + 1, 0, 0), \\
      L((a - 2)p + p - 1, p - 1, b + 2, 0) \} & \text{if } 2 \le a \text{ and } b \le p - 3.                   \\
		\end{cases}
  \end{align*}
  The different cases (and the number of irreducible constituents) correspond to the number of possible carry patterns.
\end{corollary}
\begin{proof}
By working explicitly with monomials, it is easy to read off the highest-weight vector in each
subrepresentation arising in Theorem \ref{thm:subreps-H0} (pick the highest possible
power of $x$ arising in a monomial, then the highest possible power of $y$ for this given power of
$x$).
\end{proof}
Below, we will refer to the irreducible constituent arising as a quotient of the sub-representation containing all monomials of carry pattern $c \in C(k)$ as
the $c$-carry factor.
\subsection{Global sections over the Fermat curve}
Fix $k = ap + b$ with $0 \le a, b, < p$.
We will now decompose the global sections of $\mathcal{O}_{\overline{X}_2}(k)$
into irreducibles as a $G(\Fp)$-representation. Recall that for $p$-restricted weights $\lambda$, the $G_{\Fpb}$-representation $L(\lambda)$ remains irreducible as a $G(\Fp)$-representation.

To properly account for $\delta$-twists, we need to treat $x^{p + 1} + y^{p + 1} + z^{p + 1}$ as an
element of $H^0(\mathbb{P}^2, \mathcal{O}_{\mathbb{P}^2}(p + 1) \otimes
\delta^{-p})$\footnote{This can for example be seen by computing the weight of the corresponding partial Hasse
invariant, similar to what was done for $\FlKtwo$ in §\ref{sect:2-1-hasses}.}.
If we let $j : \overline{X}_2 \rightarrow \mathbb{P}^2$ denote the
inclusion, we obtain the following short exact sequence
\begin{equation} \label{eq:SES-p2-lb}
  0 \longrightarrow \mathcal{O}_{\mathbb{P}^2}(k - p - 1) \otimes \delta^p \longrightarrow \mathcal{O}_{\mathbb{P}^2}(k) \longrightarrow j_*\mathcal{O}_{\overline{X}_2}(k) \longrightarrow 0.
\end{equation}
Since $H^1(\mathbb{P}^2, \mathcal{O}_{\mathbb{P}^2}(k)) = 0$ for any $k \in \mathbb{Z}$, we get the
following short exact sequence by taking global sections
\begin{equation} \label{eq:SES-weyl}
  0 \longrightarrow H^0(\mathbb{P}^2, \mathcal{O}_{\mathbb{P}^2}(k - p - 1) \otimes \delta^p) \longrightarrow H^0(\mathbb{P}^2, \mathcal{O}_{\mathbb{P}^2}(k)) \longrightarrow H^0(\overline{X}_2, \mathcal{O}_{\overline{X}_2}(k)) \longrightarrow 0.
\end{equation}
The following Lemma is now immediate.
\begin{lemma} \label{lem:a-eq-0}
  For $0 \le k \le p - 1$, we have
 $H^0(\overline{X}_2, \mathcal{O}_{\overline{X}_2}(k)) = L(k, 0, 0, 0)$
 as a $G(\Fp)$-representation.
\end{lemma}
The following is the main result of this subsection.
\begin{theorem} \label{thm:rep-decomposition}
	For $k = ap + b \ge p$, with $0 \le a, b < p$. We have the following decomposition in the Grothendieck
	group of $G(\Fp)$-representations.
	\[
    H^0(\overline{X}_2, \mathcal{O}_{\overline{X}_2}(k)) =
		\begin{cases}
			\begin{aligned}
        W(b, 0, -a, a) + L(p - 1, b + 1, -(a - 1), a - 1) \, + \\
        L(p - 1, p - (a - 1), b + 2, a - 2)
			\end{aligned}                                                & \text{if } a + b \le p - 1, \\
			\begin{aligned}
        W(b, 0 , -a, a) + W(p - 1, b + 1, -(a - 1), a - 1) \, + \\ L(b, p - a, 1, a - 1) - L(p - 1, b + 1, p - (a - 2), a - 2)
			\end{aligned} & \text{if } a + b \ge p.
		\end{cases}
	\]
	In particular, if $a = 1$ and $b \le p - 2$, we obtain
	\[
		H^0(\overline{X}_2, \mathcal{O}_{\overline{X}_2}(p + b)) = W(b, 0, - 1, 1) + L(p - 1, b + 1, 0, 0).
	\]
\end{theorem}
\begin{remark}
  Only the $a = 1, b \le p - 2$ case of this theorem will be used to produce
  congruences from irregular weights.
  We include this more general result, as this was useful in earlier stages of the project, and we
  believe that it could be of independent interest.
\end{remark}
\begin{proof}[Proof of theorem]
  From \eqref{eq:SES-weyl}, we get the following equality in the Grothendieck group
  of $G(\Fp)$-representations (using
  that $\delta^p = \delta$ as a character of $G(\Fp)$).
	\begin{equation}
    H^0(\overline{X}_2, \mathcal{O}_{\overline{X}_2}(ap + b)) = W(ap + b, 0, 0, 0) - W((a - 1)p + (b - 1), 0, 0, 1).
		\label{eq:total-difference}
	\end{equation}
  Corollary \ref{cor:weyl-decomp} implies that this difference decomposes further,
  giving up to three positive, and up to three negative terms (depending on $|C(k)|$ and $|C(k - (p
  + 1))|$).
  To make calculation easier, we will consider the differences between corresponding carry-factors appearing in 
  \eqref{eq:total-difference} separately. These individual differences will a priori only
	make sense in the Grothendieck group, and it is only their total sum that gives an honest
	representation in general. Furthermore, it is entirely possible that the positive and negative term of \eqref{eq:total-difference} above might
	not admit carry-factors of the same carry-pattern simultaneously, so we will need to take this into
	account during the calculations below.
	
	Let us start with the terms of \eqref{eq:total-difference} corresponding to $0$-carrying,
  these terms will always exist for $k \ge p$ (but note that for $b = 0$ the negative term
  vanishes).
  Hence, we are lead to consider the difference
  \[
    L(ap + b, 0 , 0, 0) - L((a - 1)p + (b - 1), 0, 0, 1).
  \]
  Using the Steinberg tensor product theorem \cite[II.3.17]{jantzen}, we can write this as
	\begin{equation} \label{eq:0-carry-diff}
    L(b, 0, 0, 0) \otimes L(0, 0, -a, a) - L(b - 1, 0, 0) \otimes L(0, 0, -(a - 1), a).
	\end{equation}
  Let us first decompose the positive term. Since $a < p$, we have that $L(0, 0, -a, a) \cong W(0, 0, -a, a)$ and hence we can apply Klimyk's formula for decomposing the tensor product \cite[II.5.8b]{jantzen}. To apply it, we need to know the weight multiplicities of $L(b, 0, 0, 0)$. Since $b < p$, we can describe
  this representation as the space of degree $b$ homogeneous polynomials. In particular, we can read
  off the weight multiplicities from the monomials in this space.
	For a given weight $\lambda \in X(T)$, we let $\chi(\lambda)$ denote the character of $W(\lambda)$,
  and let $e(\lambda)$ denote the character of the 1-dimensional representation $\lambda$. Thus, applying Klimyk's formula gives us (where the indices of the sums appearing are always non-negative)
	\begin{align*}
		W(0, 0, -a, a) \otimes L(b, 0, 0, 0) & = \chi(0, 0, -a, a) \sum_{x + y + z = b} e(x, y, z, 0) \\
		                               & = \sum_{x + y + z = b} \chi(x, y, z - a, a)         \\
		                               & = \sum_{x + y \le b} \chi(x, y, b - x - y - a, a).
	\end{align*}
	Applying the same procedure to the negative term, we obtain
	\begin{align*}
		W(0, 0, -(a - 1), a) \otimes L(b - 1, 0, 0, 0) & = \sum_{x + y \le b - 1} \chi(x, y, b - 1 - x - y - (a - 1), a) \\
		                                         & = \sum_{x + y \le b - 1} \chi(x, y, b - x - y - a, a).
	\end{align*}
	The difference between these two sums is then
	\[
		\sum_{x + y = b} \chi(x, y, b - x - y - a, a).
	\]
	To simplify this sum further, note that for any character $\mu$ and any root
  $\alpha$ we have $\chi(s_\alpha \cdot \mu) = - \chi(\mu)$ (where the action here is the dot action). Note that $s_{\alpha_1} \cdot (x, y, b - x - y - a, a) = (y - 1, x + 1, b - x - y - a, a)$. This implies that, if $x$ is not maximal (i.e. $x \neq b$), then there exists a unique triple $(y' - 1, x' + 1, z', a)$ in the sum that cancels this term. Hence, after cancellations, only the term with $x = b$ remains, and we see that the
	difference in \eqref{eq:0-carry-diff} becomes
	\[
		L(0, 0, -a, a) \otimes L(b, 0, 0, 0) - L(0, 0, -(a - 1), a) \otimes L(b - 1, 0, 0, 0) = W(b, 0, -a, a).
	\]
  In particular this is always an honest representation, and has either one or two irreducible constituents, depending on whether $(b, 0, -a, a)$ lies in the lower or upper alcove (i.e. depending on if $a + b \le p - 2)$ or not). This gives us the first term in both cases of the theorem.
	
	Let us now apply the same procedure to the terms of \eqref{eq:total-difference} corresponding to $1$-carrying.
	In this case, note that $W((a - 1)p + (b - 1), 0, 0, 0)$ does not admit $1$-carrying if $a = 1$. Thus, for $a \ge 2$, we are interested in studying the difference 
\[
  L((a - 1)p + p - 1, b + 1, 0, 0) - L((a - 2)p + p - 1, b, 0, 1).
\]
Applying the Steinberg tensor product theorem gives
\begin{equation} \label{eq:1-carry-diff}
  L(p - 1, b + 1, 0, 0) \otimes L(0, 0, -(a - 1), a - 1) - L(p - 1, b, 0, 0) \otimes L(0, 0, -(a - 2), a - 1).
\end{equation}
However, note that we can use this expression to treat the $a \ge 2$ and $a = 1$ cases simultaneously,
as the negative term vanishes for $a = 1$. Thus, we can proceed without case-splits.

We start by considering the first term of \eqref{eq:1-carry-diff}. We have that $L(p - 1, b + 1, 0, 0)$ corresponds to the quotient of $W(p + b, 0, 0, 0)$ spanned by the monomials such that the sum of their exponents in base $p$ has 1-carrying. Hence, its weight spaces all have multiplicity 1, and are given by the non-negative tuples $(x, y, z, 0)$ such that $x + y + z = p + b$ and $x, y, z \le p - 1$. This gives us the following equality in the Grothendieck group
	\begin{align*}
		L(p - 1, b + 1, 0, 0) \otimes L(0, 0, -(a - 1), a - 1) & = \sum_{\substack{x, y, z \le p - 1 \\ x + y + z = p + b}} \chi(x, y, z - (a - 1), a - 1) \\
		                                             & = \sum_{\substack{x, y \le p - 1    \\ b + 1 \le x + y \le p + b}} \chi(x, y, p + b - x - y - (a - 1), a - 1).
	\end{align*}
	Similarly, we obtain
	\[
		L(p - 1, b, 0, 0) \otimes L(0, 0, -(a - 2), a - 1) = \sum_{\substack{x, y \le p - 1 \\ b \le x + y \le p + b - 1}} \chi(x, y, p + b - x - y - (a - 1), a - 1).
	\]
	Hence, the difference is given by
	\[
		\sum_{\substack{x, y \le p - 1 \\ x + y = p + b}} \chi(x, y, -(a - 1), a - 1) - \sum_{\substack{x, y \le p - 1 \\ x + y = b}} \chi(x, y, p - (a - 1), a - 1).
	\]
	As before, we get cancellation of all the terms in the two sums, except when $x$ takes on the maximal value $p - 1$ (resp. $b$). Thus, the expression simplifies to
	\[
		\chi(p-1, b + 1, -(a - 1), a - 1) - \chi(b, 0, p - (a - 1), a - 1).
	\]
  Applying the dot action with respect to $s_{\alpha_2}$ to the negative term gives us
	\[
		\chi(p - 1, b + 1, -(a - 1), a - 1) + \chi(b, p - a, 1, a - 1).
	\]
  If $a + b > p$ then the latter term lies in the lowest alcove, but the first term is not $p$-restricted, and we do not know its decomposition factors as a $G(\Fp)$-representation. This gives us the second and third terms in the case $a + b \ge p$ of the theorem. Otherwise, if $a + b \le p - 1$ then we can apply the dot action of the simple reflection $s_{\alpha_1}$ and obtain
	\[
		\chi(p - 1, b + 1, -(a - 1), a - 1) - \chi(p - a - 1, b + 1, 1, a - 1).
	\]
  Notice that the character of the positive term is obtained by reflecting the character of the negative term from the lower to the upper alcove. Hence, in this case, we get that the difference is equal to $L(p - 1, b + 1, -(a - 1), a - 1)$ in the Grothendieck group. This gives us the second term of the theorem in the case $a + b \le p - 1$.
	
	If $a = 1$, we have now proved the theorem, as there is no $2$-carrying. So let us proceed by assuming $a \ge 2$, and consider
	the terms of \eqref{eq:total-difference} corresponding to 2-carrying.
  In other words, we wish to study (ignoring edge cases for now)
  \[
    L((a - 2)p + p - 1, p - 1, b + 2, 0) - L((a - 3)p + p - 1, p - 1, b + 1, 1).
  \]
  Applying the Steinberg tensor product theorem gives
	\[
    L(p - 1, p - 1, b + 2, 0) \otimes L(0, 0, -(a - 2), a - 2) - L(p - 1, p - 1, b + 1, 0) \otimes L(0, 0, -(a - 3), a - 2).
	\]
	As before, notice that the edge cases $a = 2$ and $b = p - 2$ can be treated by the same formula, as the
	corresponding terms vanish in these cases. Hence, we can proceed without case splits.
	
	Let us start by studying the first term. Recall that $L(p - 1, p - 1, b + 2, 0)$ corresponds to the quotient of $W(2p + b, 0, 0, 0)$ spanned by the monomials such that the sum of their exponents in base $p$ has $2$-carrying. Hence, its weight spaces all have multiplicity 1, and the characters are given by tuples $(x, y, z, 0)$ such that $x + y + z = 2p + b$ and $x, y, z \le p - 1$. As before, we get
	\begin{align*}
		L(p - 1, p - 1, b + 2, 0) \otimes L(0, 0, -(a - 2), a - 2) & = \sum_{\substack{x, y, z \le p - 1 \\ x + y + z = 2p + b}} \chi(x, y, z - (a - 2), a - 2) \\
		                                                 & = \sum_{\substack{x, y \le p - 1    \\ p + b + 1 \le x + y}} \chi(x, y, 2p + b - x - y - (a - 2), a - 2).
	\end{align*}
	Similarly, we get for the negative term
	\[
		L(p - 1, p - 1, b + 1, 0) \otimes L(0, 0, -(a - 3), a - 2) = \sum_{\substack{x, y \le p - 1 \\ p + b \le x + y}} \chi(x, y, 2p + b - x - y - (a - 2), a - 2).
	\]
	Thus, the difference is given by
	\[
		- \sum_{\substack{x, y \le p - 1 \\ x + y = p + b}} \chi(x, y, p - (a - 2), a - 2).
	\]
  As before, we get cancellation of all terms where $x$ is not maximal, i.e. $x \neq p - 1$. So the only term that remains is $- \chi(p - 1, b + 1, p - (a - 2), a - 2)$. The character here is dominant whenever $a + b \ge p$, giving us only a virtual representation. This concludes the $a + b \ge p$ case of the theorem. Otherwise, if $a + b < p$, applying the dot action of $s_{\alpha_2}$ gives us
	\[
		- \chi(p - 1, b + 1, p - (a - 2), a - 2) = \chi(p - 1, p - (a - 1), b + 2, a - 2).
	\]
  This proves the remaining cases of the decomposition.
\end{proof}

\subsection{Proof of Proposition \ref{prop:main-decomp}} \label{sect:global-sect-zcb}
Now we apply the results of the previous section to prove Proposition \ref{prop:main-decomp}.
\begin{proof}[Proof of Proposition \ref{prop:main-decomp}]
  Our goal is to compute the global sections of $\mathcal{O}_{\overline{X}_1}(\lambda')$.
  Via the isomorphism $G / P \cong G / P^*$ (sending $g \mapsto w_0g^{-t}w_0$), this corresponds to the following line bundle, where
  we have also factored out determinant and multiplier twists
  \begin{equation} \label{eq:twisted}
  \mathcal{O}_{\overline{X}_2}(2p - (a - b + 2)) \otimes {\det}^{\otimes - (b - 2)} \otimes
    \delta^{-(2p + d)}.
  \end{equation}
  Ignoring the twists, we can apply Theorem \ref{thm:rep-decomposition} when $a - b \le p -
  2$ to get the exact sequence
  \[
    0 \longrightarrow W(p - (a - b + 2), 0, -1, 1) \longrightarrow
    H^0(\overline{X}_2, \mathcal{O}_{\overline{X}_2}(2p - (a - b + 2)))
    \longrightarrow L(p - 1, p - 1 - (a - b), 0, 0) \longrightarrow 0.
  \]
  The fact that $W(p - (a - b + 2), 0, -1, 1)$ is actually a subrepresentation can be seen from the proof of Theorem \ref{thm:rep-decomposition}, and the fact that, for $a - b \le p - 2$, we have that the morphism $W((a - 1)p + (b - 1), 0, 0, 1) \to W(ap +
  b, 0, 0, 0)$ given by multiplying by $x^{p + 1} + y^{p + 1} + z^{p + 1}$ can only decrease the carry-pattern.
  The proposition now follows in the $a - b \le p - 2$ case by incorporating the determinant and
  multiplier twists from \eqref{eq:twisted}, switching back from $\overline{X}_2$ to $\overline{X}_1$ (which has the effect $(a, b, c, d) \mapsto (-c, -b, -a, -d)$ on weights), and finally shifting by the trivial weight $(p + 1, p + 1, p + 1, 0)$ (corresponding to the character $\det^{p + 1}$, which is trivial on $G(\Fp)$).
  The $a - b = p - 1$ case follows similarly from Lemma \ref{lem:a-eq-0}.
\end{proof}

\section{The exotic congruence} \label{sect:applications}
\subsection{The Hecke algebra}
Let $Q$ denote the set containing all primes $\ell$ such that either $G$ (or equivalently $G'$) ramifies at $\ell$, or such that
$K$ is not hyperspecial at $\ell$.
Consider the following global Hecke algebra, defined as a restricted tensor product
\[
  \mathbf{T} := {\bigotimes}_{\ell \not \in Q \cup \{ p \}}' \Z_{\ell}\left[ G(\Z_\ell) \backslash G(\Q_\ell) / G(\Z_\ell)  \right].
\]
A \textit{mod $p$ Hecke eigensystem} is defined to be a morphism $\mathbf{T} \rightarrow \Fpb$.
Given a $\Fpb$-vector space $M$, equipped with a $\mathbf{T}$-action,
any simultaneous eigenvector for the $\mathbf{T}$-action on $M$ gives
rise to a mod $p$ Hecke eigensystem. We will say that the mod $p$ Hecke eigensystem
\textit{appears in $M$}.
\subsubsection{Action on coherent cohomology}
The global Hecke algebra $\mathbf{T}$ acts on the coherent cohomology groups
$H^0(\Sh_{K, \Fpb}, \omega(\lambda)) = H^0(\Fl_{K, \Fpb}, \mathcal{L}(\lambda))$ via
Hecke correspondences (and similarly in the type $(3, 0)$-case). See \cite[§6.1]{goldring} for details.

\subsubsection{Restriction to strata}
The correspondences giving rise to the action of $\mathbf{T}$ involve parametrizing certain
prime-to-$p$ isogenies.
In particular, since lying on a certain stratum is a condition on the $p$-torsion,
the strata will be invariant under these correspondences.
In particular, given any coherent weight $\lambda$ and any stratum $\overline{\mathcal{Z}}_w \subseteq \Fl_K$, the space $H^0(\overline{\mathcal{Z}}_{w, \Fpb}, \mathcal{L}(\lambda))$ is a $\mathbf{T}$-module, and the corresponding restriction maps are Hecke-equivariant.

\subsection{Injectivity of restriction}
Let us now show injectivity on global sections of the restriction map to $\Zcb_{2}$ for irregular
weights.
\begin{proposition}
  For $p > 3$, the map $H^0(\Fl_{K, \Fpb}, \mathcal{L}(\lambda)) \rightarrow H^0(\Zcb_{2, \Fpb}, \mathcal{L}(\lambda))$ is
  injective for all irregular weights $\lambda = (a, b, b - 2, d)$ with $0 \le a - b \le p - 1$.
  For $p = 3$, injectivity holds for $a - b = p - 1 = 2$.
\end{proposition}
\begin{proof}
  Since twisting by the multiplier character does not affect injectivity, we will
  for simplicity ignore the multiplier factor and treat weights as elements of $\Z^{\oplus 3}$. Furthermore, up to twisting by
  determinants, we may assume that $b = 0$, so we can for simplicity work with the irregular weight
  $\lambda = (k, 0, -2)$ with $0 \le k \le p - 1$.

  Let us first show that the restriction map to $\Zcb_{12}$ is injective.
  Let $i_{12} : \Zcb_{12} \to \Fl_K$ denote the inclusion.
  By Lemma \ref{lem:vanish-simply}, the partial Hasse invariant $\ha_{12} \in H^0(\Fl_K, \mathcal{L}(p, 0, 1))$ vanishes simply on $\Zcb_{12}$. We have the following exact sequence
  \[
    0 \longrightarrow \mathcal{L}(\lambda - (p, 0, 1)) \longrightarrow \mathcal{L}(\lambda)
    \longrightarrow i_{12, *}(\mathcal{L}(\lambda)|_{\Zcb_{12}}) \longrightarrow 0.
  \]
  Thus, it remains to show that $H^0(\Fl_{K, \Fpb}, \mathcal{L}(k - p, 0, -3))$ is zero.
  Notice that along any fiber of $\FlKtwo
  \rightarrow \ShKtwo$ (any of which is isomorphic to $\mathbb{P}^1$), we have that
  $\mathcal{L}(k - p, 0, -3)|_{\mathbb{P}^1} \cong \mathcal{O}_{\mathbb{P}^1}(k - p)$. But $k - p < 0$, so
  the global sections must vanish.
  
  It remains to show that $H^0(\Zcb_{12, \Fpb}, \mathcal{L}(\lambda)) \rightarrow H^0(\Zcb_{2, \Fpb}, \mathcal{L}(\lambda))$ is injective.
  Assume that $f \in H^0(\Zcb_{12, \Fpb}, \mathcal{L}(\lambda))$ vanishes along $\Zcb_{2, \Fpb}$.
  Recall from Lemma \ref{lem:P1-fibers} that $\Zcb_{1}$ is a disjoint union of $\mathbb{P}^1$s.
  The proof of \cite[Proposition 2.10]{gorenshalit},
  shows that
  $\Zcb_{2}$ meets each $\mathbb{P}^1$-fiber of $\Zcb_{1}$ at $p + 1$ distinct points.
  We know that $\mathcal{L}(\lambda)|_{\Zcb_{1}}$ is isomorphic to $\mathcal{O}_{\mathbb{P}^1}(k)$ along each fiber.
  Hence, since $f$ vanishes on $\Zcb_{2}$, it must vanish at $p + 1$ points of each $\mathbb{P}^1$.
  Since $k < p + 1$, this implies that $f$ must also vanish on $\Zcb_{1}$.
  Thus, by Lemma \ref{lem:vanish-simply}, $f$ is divisible by $\ha_{2} \cdot \ha_{1} \in H^0(\Zcb_{2, \Fpb}, \mathcal{L}(p + 1, p, 1))$.
  Thus, after dividing by $\ha_{2} \cdot \ha_{1}$, it suffices to show that any section in $H^0(\Zcb_{12, \Fpb}, \mathcal{L}(\lambda - (p + 1, p, 1))$ must vanish.
  Note that $\lambda - (p + 1, p, 1) = (k - p - 1, - p, - 3)$.
  One easily observes that the degree along the $\mathbb{P}^1$-fibers is now $k - 1$.
  Thus, dividing by $\ha_{2} \cdot \ha_{1}$ decreases the degree along the $\mathbb{P}^1$-fibers.
  Next, it follows from Theorem \ref{thm:aut-vb-descr}
  that (up to twisting) we have $\mathcal{L}((k - p - 1, - p, - 3))|_{\Zcb_{2}} \cong \mathcal{L}'_{P^*}(0, 0, -((3 - p)p - (k + 3 - p - 1)))|_{\Zcb'_{1}}$.
  Recall that $\Zcb'_{1}$ is just a disjoint union of Fermat curves.
  One calculates that the degree of the corresponding line bundle along each Fermat curve is $(p + 1)(-p^2 + 4p - (k + 2))$.
  For $p > 3$ (or $k = p - 1$ for $p = 3$), this has negative degree.
  Since the degree is now negative, any section of $\mathcal{L}((k - p - 1, - p, -
  3))|_{\Zcb_{12}}$, must vanish along $\Zcb_{2}$. By the same argument as before,
  this implies that it must also vanish along $\Zcb_{1}$. So any section of $\mathcal{L}((k - p - 1, - p, -
  3))|_{\Zcb_{12}}$ is divisible by $\ha_{2} \cdot \ha_{1}$.
  However, note that we have also shown that
  division by $\ha_{2} \cdot \ha_{1}$ decreases the degree along both
  $1$-dimensional strata, so we can repeat this argument indefinitely.
  This implies that any section of $\mathcal{L}((k - p - 1, - p, -
  3))|_{\Zcb_{12}}$ must vanish to infinite degree along $V(\ha_{2}\cdot \ha_{1})$,
  so it must be zero.
  Thus, we are done.
\end{proof}
In particular, we obtain the following immediate corollary.
\begin{corollary} \label{cor:injectivity}
  Fix an irregular weight $\lambda = (a, b, b - 2, d)$ with $0 \le a - b \le p - 1$, if $p = 3$ we further
  assume that $a - b = p - 1$.
  Then every mod $p$ Hecke eigensystem appearing in $H^0(\Fl_{K, \Fpb}, \mathcal{L}(\lambda))$
  also appears in $H^0(\Zcb_{2, \Fpb}, \mathcal{L}(\lambda))$.
\end{corollary}

\subsection{The congruence}
Let us now show the main theorem of this paper.
\begin{theorem} \label{thm:main-thm}
  Fix an irregular weight $\lambda = (a, b, b - 2, d)$ with $0 \le a - b \le p - 1$.
  Set $\lambda_1 = (p + b - 1, a, b, d + 2)$ and $\lambda_2 = (p + b, p + b - 1, a + 1, d + 1)$.
  If $a - b \le p - 2$ and $p > 3$, then every mod $p$ Hecke eigensystem appearing in $H^0(\Fl_{K, \Fpb}, \mathcal{L}(\lambda))$, must appear in
  either $H^0(\Fl'_{K', \Fpb}, \omega'(L(\lambda_1)))$ or
  $H^0(\Fl'_{K', \Fpb}, \omega'(W(\lambda_2)))$.
  If $a - b = p - 1$, then the corresponding mod $p$ Hecke eigensystem must appear
  in $H^0(\Fl'_{K', \Fpb}, \omega'(L(\lambda_1)))$.
\end{theorem}
\begin{proof}
  This follows from Corollary \ref{cor:injectivity} and Theorem
  \ref{thm:main-decomp}, provided we know
  that the identification $\Zcb_{2} \cong \Zcb'_1$ respects the corresponding
  Hecke operators on both sides.
  But this is clear, since the identification $\Zcb_{2} \cong \Zcb'_1$ is given by a $p$-power isogeny on the corresponding universal abelian varieties, so in particular
  it does not affect the action of $\mathbf{T}$, which only involves isogenies of degree coprime to
  $p$.
\end{proof}

\printbibliography

\end{document}